\documentclass{article}
\usepackage{arxiv}
\usepackage{mathptmx}
\DeclareMathAlphabet{\mathcal}{OMS}{cmsy}{m}{n}
\RequirePackage{fix-cm}
\usepackage[sort&compress]{natbib}
\usepackage{siunitx}
\usepackage{amsmath,bm,mathtools}

\usepackage[utf8]{inputenc}
\usepackage[T1]{fontenc}
\usepackage[outline]{contour}
\contourlength{1.2pt} 
\usepackage{graphicx}
\usepackage[figurename=Fig.,labelfont=bf,labelsep=period]{caption}
\usepackage{subcaption}
\usepackage{amsfonts}
\usepackage{amssymb}
\usepackage{amsbsy}
\usepackage[percent]{overpic}
\usepackage[colorlinks=true,citecolor=red,linkcolor=black]{hyperref}

\usepackage{url}            
\usepackage{booktabs}       
\usepackage{nicefrac}       
\usepackage{microtype}      
\usepackage{lipsum}
\usepackage{algorithm,tikz}
\usepackage[]{algpseudocode}
\usepackage{pgfplots}
\usepackage{subcaption}
\usepackage[]{threeparttable}
\usepackage{multirow,bm}
\usepackage{gensymb}

\newcommand{\Am}{\mathbf{A}} 
\newcommand{\am}{\mathbf{a}}
\newcommand{\Bm}{\mathbf{B}}

\newcommand{\Dm}{\mathbf{D}}

\newcommand{\Dtr}{\mathcal{D}_\mathrm{tr}}
\newcommand{\Dval}{\mathcal{D}_\mathrm{val}}

\newcommand{\Gm}{\mathbf{G}}
\newcommand{\Hb}{\mathbf{H}_\mathbf{B}}
\newcommand{\Hl}{\mathbf{H}_\mathbf{L}}

\newcommand{\Hm}{\mathbf{H}}
\newcommand{\Lm}{\mathbf{L}}

\newcommand{\Mm}{\mathbf{M}_\mathrm{s}}

\newcommand{\Ntr}{{N}_\mathrm{tr}}
\newcommand{\Nval}{{N}_\mathrm{val}}
\newcommand{\ngm}{n_\mathrm{g}}

\newcommand{\ny}{n_\mathrm{y}}

\newcommand{\nw}{n_\mathrm{w}}
\newcommand{\no}{n_\mathrm{o}}
\newcommand{\ppm}{\mathbf{p}}
\newcommand{\pnn}{\mathbf{p}_{\!_\mathrm{NN}}}
\newcommand{\Qy}{Q_\mathrm{y}}
\newcommand{\rr}{\mathbf{r}}

\newcommand{\Km}{\mathbf{K}_\mathrm{s}}

\newcommand{\kpost}{k_\mathrm{post}}
\newcommand{\kpre}{k_\mathrm{pre}}
\newcommand{\Cm}{\mathbf{C}}
\newcommand{\Cms}{\mathbf{C}_\mathrm{s}}

\newcommand{\gm}{\mathbf{g}}

\newcommand{\Sm}{\mathbf{S}}

\newcommand{\wm}{\mathbf{w}}
\newcommand{\um}{\mathbf{u}}

\newcommand{\Xm}{\mathbf{X}}
\newcommand{\Xmb}{\overline{\mathbf{X}}}
\newcommand{\xmb}{\overline{\mathbf{x}}}
\newcommand{\Ym}{\mathbf{Y}}
\newcommand{\xm}{\mathbf{x}}
\newcommand{\xo}{\overline{\xm}}
\newcommand{\xnl}{\mathbf{x}_\mathrm{corr}}

\newcommand{\ym}{\mathbf{y}}

\newcommand{\drm}{\mathrm{d}}
\newcommand{\transpose}{^{\mathrm{T}}}

\newcommand{\xii}{\pmb{\xi}}

\newcommand{\thetaa}{\pmb{\theta}}

\newcommand{\vey}{\mathrm{vec}(\mathcal{Y})} 
\newcommand{\veyn}{\mathrm{vec}(\mathcal{Y}_\mathrm{NN})}

\newcommand{\vecy}{\mathrm{vec}\left(\mathcal{Y}(\thetaa,\xii)\right)}

\newcommand{\Exp}{\mathbb{E}}

\newcommand\Tstrut{\rule{0pt}{2.6ex}}         
\newcommand\Bstrut{\rule[-0.9ex]{0pt}{0pt}}   

\newtheorem{theorem}{Proposition}

\usepackage{tikz}
\usetikzlibrary{matrix,chains,positioning,decorations.pathreplacing,arrows}
\usetikzlibrary{positioning,calc}
\usetikzlibrary{patterns,decorations.pathmorphing,decorations.markings}
\usetikzlibrary{shapes}
\usetikzlibrary{fit}
\tikzstyle{block} = [draw,rectangle,thick,minimum height=2em,minimum width=2em]
\tikzstyle{sum} = [draw,circle,inner sep=0mm,minimum size=2mm]
\tikzstyle{connector} = [->,thick]
\tikzstyle{line} = [thick]
\tikzstyle{branch} = [circle,inner sep=0pt,minimum size=1mm,fill=black,draw=black]
\tikzstyle{guide} = []
\tikzset{%
	every neuron/.style={
		circle,
		draw,
		fill = green!20,
		thick,
		minimum size=1.1cm
	},
	neuron missing/.style={
		draw=none, 
		scale=4,
		fill = none,
		text height=0.333cm,
		execute at begin node=\color{black}$\vdots$
	},
	input neuron/.style={
		circle,
		draw,
		thick,
		fill = yellow!20,
		minimum size=1cm
	},
	output neuron/.style={
		circle,
		draw,
		thick,
		fill = red!20,
		minimum size=1cm
	},
}

\usetikzlibrary{3d,decorations.text,shapes.arrows,positioning,fit,backgrounds}
\tikzset{pics/fake box/.style args={
		#1 with dimensions #2 and #3 and #4}{
		code={
			\draw[gray,ultra thin,fill=#1]  (0,0,0) coordinate(-front-bottom-left) to
			++ (0,#3,0) coordinate(-front-top-right) --++
			(#2,0,0) coordinate(-front-top-right) --++ (0,-#3,0) 
			coordinate(-front-bottom-right) -- cycle;
			\draw[gray,ultra thin,fill=#1] (0,#3,0)  --++ 
			(0,0,#4) coordinate(-back-top-left) --++ (#2,0,0) 
			coordinate(-back-top-right) --++ (0,0,-#4)  -- cycle;
			\draw[gray,ultra thin,fill=#1!80!black] (#2,0,0) --++ (0,0,#4) coordinate(-back-bottom-right)
			--++ (0,#3,0) --++ (0,0,-#4) -- cycle;
			\path[gray,decorate,decoration={text effects along path,text={}}] (#2/2,{2+(#3-2)/2},0) -- (#2/2,0,0);
		}
}}
\tikzset{circle dotted/.style={dash pattern=on .05mm off 2mm,
		line cap=round}}

\title{Bi-fidelity Design of Nonlinear Dynamical Systems Under Uncertainty Using Neural Networks} 

\author{
Maryam Maghazeh\\
  Department of Mechanical Engineering\\
  Northern Arizona University\\
  Flagstaff, AZ 86011 \\
  \texttt{mm5335@nau.edu} \\
\And 
 {Patrick T. Brewick} \\
	Department of Civil and Environmental\\  Engineering and Earth Sciences\\
	University of Notre Dame\\
	Notre Dame, IN 46556 \\
	\texttt{pbrewick@nd.edu} \\ 
    \And
    Subhayan De\footnote{Corresponding author}\\
  Department of Mechanical Engineering\\
  Northern Arizona University\\
  Flagstaff, AZ 86011 \\
  \texttt{Subhayan.De@nau.edu} \\
}

\begin{document}

\maketitle

\begin{abstract} 

Design of nonlinear dynamical systems under uncertainty is computationally expensive since each optimization iteration often requires thousands of transient simulations and sensitivity evaluations. This paper develops a bi-fidelity approach for design under uncertainty of dynamical systems with nonlinear behavior in spatially localized components within an otherwise linear model. The system response is separated into a nominal linear contribution and a correction generated by an equivalent pseudoforce that represents the localized nonlinear and uncertain behavior. A neural network is trained to approximate this pseudoforce, while the known linear state-space model propagates its effect through the full system. Design sensitivities are obtained by differentiating the neural network and passing the resulting derivatives through the linear dynamics. Theoretical bounds are also  developed to relate errors in the learned pseudoforce and its design derivatives obtained through neural networks to errors in the response, objective, constraints, and their design gradients. The method is demonstrated on a nonlinear solar heating system and a base-isolated structure with hysteretic dynamics. 
In both examples, the proposed approach accurately approximates the pseudoforce, its design derivative on validation cases, and  objective and constraint statistics at the optimized designs. 
Once the offline training is completed, it substantially reduces the online cost of gradient-based design under uncertainty. 
\end{abstract}

\section{Introduction}

Nonlinear dynamical models are commonly used in structural, mechanical, thermal, and energy systems whenever restoring forces, dissipation, transport coefficients, or constitutive behavior depend on the evolving state. 
As a result, their performance cannot usually be described by a single equilibrium state. Quantities used in design of these systems, such as peak displacement, vibration energy, maximum temperature, or accumulated thermal output, therefore depend on the transient trajectory of the system \cite{nayfeh2024nonlinear,nayfeh2008applied,strogatz2024nonlinear,kerschen2006past}. Many of these nonlinear systems are high-dimensional but contain nonlinearities in only a small number of components. A base-isolated building, for example, may have a large linear superstructure coupled to a low-dimensional hysteretic isolation layer \cite{wen1976method,ramallo2002smart,kamalzare2015efficient}. Similar configurations arise in mechanical assemblies with nonlinear joints \cite{ferri1988modeling,bowden1990joint}, nonlinear inerters \cite{li2025model}, and in structures equipped with localized passive control devices \cite{balaji2021applications,de2017efficient,sarnaik2025design}. Although the nonlinear behavior is localized, its effect propagates through the entire system. 

Designing such systems requires repeated evaluation of their transient response as the design parameters change. For example, the objective may depend on an integrated quantity, such as vibration energy or thermal output, while the constraints may involve peak displacement, maximum temperature, or another extreme response over the simulation interval \cite{nocedal2006numerical,haftka2012elements,arora2004introduction}. Gradient-based optimization also requires the sensitivity of these quantities with respect to the design variables. For nonlinear dynamical systems, these sensitivities can be evaluated using finite differences, direct differentiation, adjoint formulations, or algorithmic differentiation \cite{cao2003adjoint,griewank2008evaluating,wojtkiewicz2014efficient}. The cost is manageable for a small deterministic model, but it grows quickly for large systems and long time histories.

Additionally, uncertainty in material properties, damping, nonlinear constitutive parameters, loading, and environmental conditions further complicates the design of these nonlinear dynamical systems. Even modest parameter variations can change the amplitude, phase, or duration of a nonlinear response. Design under uncertainty must therefore account for a collection of trajectories rather than a single nominal simulation \cite{ghanem2003stochastic,le2010spectral,smith2024uncertainty,sullivan2015introduction}. Robust design commonly uses statistics such as the mean and variance of a response quantity \cite{beyer2007robust,de2020topology}, whereas reliability-based design constrains the probability of exceeding a prescribed limit \cite{du2004sequential,hu2024reliability,maghazeh2026reliability}. Design under uncertainty has also been considered for coupled structural and vibration-control systems with interval-valued parameters \cite{yang2025integrated}. However, all of these require repeated evaluations of the governing equations as the design changes. 
The resulting cost can become prohibitive when the system is high-dimensional, since the response must be resolved over long time intervals, or when a large number of uncertainty realizations is needed at each iteration.


Surrogate models offer one way to reduce this computational burden by replacing repeated nonlinear simulations with inexpensive approximations \cite{jiang2019surrogate,forrester2007multi}. 
Among these, data-driven models have been used extensively to approximate nonlinear dynamics. Representative approaches include residual neural networks \cite{de2022prediction}, recurrent neural networks \cite{salmela2021predicting,li2021recurrent}, neural ordinary differential equations \cite{chen2018neural}, neural network-based reduced-order models and autoencoders for transient nonlinear evolution \cite{san2019artificial,champion2019data,xu2026pod}, Koopman-based embeddings and physics-integrated formulations incorporating nonlinear normal modes \cite{lusch2018deep,rostamijavanani2023study}, sparse and physics-informed methods for discovering governing equations \cite{brunton2016discovering,saqlain2023discovering}, Gaussian
process-based dynamical models \cite{fezai2026data,ye2024gaussian}, and neural operators \cite{lu2021learning,li2020fourier}. 
In multi-fidelity modeling, the cost of uncertainty propagation and optimization can be reduced using an inexpensive approximate model combined with a smaller set of high-fidelity data   \cite{kennedy2000predicting,peherstorfer2018survey}. Recently, neural networks have been used to learn nonlinear corrections between fidelity levels for parametric models and differential equations \cite{perdikaris2017nonlinear,meng2020composite,de2020transfer}. 

Another avenue for reducing the cost of nonlinear dynamical simulations is to decompose the response 
into the trajectory of a nominal linear model and a correction induced by an equivalent force representing the nonlinear and uncertain terms \cite{wojtkiewicz2011efficient,wojtkiewicz2014efficient,de2018computationally}. This equivalent force, referred to herein as the nonlinear pseudoforce, is typically much lower-dimensional than the full state. However, its exact evaluation leads to a nonlinear Volterra integral equation of the second kind in a nonstandard form and requires an iterative solution for each combination of design and uncertain parameters \cite{wojtkiewicz2011efficient,linz1985analytical,brunner2004collocation,de2018computationally}. 
Volterra-series representations have been used for nonlinear identification of multi degree of freedom dynamical systems \cite{prawin2017nonlinear}. 
Neural networks have also been used to identify nonlinear input-output representations, including the recovery of multi-input multi-output Volterra kernels in the frequency domain \cite{preston2026mimo}. 
In a previous work, De et al. \cite{de2020uncertainty} used a neural network to approximate the pseudoforce and combined it with the nominal linear response for uncertainty propagation. That study addressed uncertainty propagation at prescribed system configurations. It did not, however, consider a pseudoforce representation over a design domain, automatic differentiation of the learned pseudoforce with respect to design variables, propagation of these derivatives to optimization gradients, or gradient-based design under uncertainty. These issues constitute the focus of the present work.


Using the pseudoforce decomposition, the present work constructs a physics-guided bi-fidelity framework for the design of nonlinear dynamical systems under uncertainty. The two fidelity levels are coupled through the equivalent nonlinear forcing rather than through direct interpolation of their response predictions. The nominal linear system provides the low-fidelity model for system-level propagation, while high-fidelity nonlinear simulations supply the pseudoforce data used to train a neural network over the joint design and uncertainty domains. The learned forcing is combined with the nominal dynamics to reconstruct nonlinear transient responses. Its derivatives with respect to the design variables are obtained by automatic differentiation and propagated through the same nominal dynamics to calculate response sensitivities \cite{wojtkiewicz2011efficient,wojtkiewicz2014efficient,paszke2017automatic,griewank2008evaluating}. The work makes four principal contributions:
\begin{itemize}
    \item It constructs a physics-guided bi-fidelity representation in which the nominal linear model retains the known system dynamics, while high-fidelity simulations inform a learned correction for the localized nonlinear behavior.
    \item It develops a design- and uncertainty-dependent representation of the nonlinear pseudoforce that is not tied to a particular response output or design functional.
    \item It establishes a differentiable procedure for reconstructing nonlinear transient responses and their design sensitivities without repeatedly solving the full nonlinear governing equations during uncertainty propagation and optimization.
    \item It derives error bounds that characterize how inaccuracies in both the learned nonlinear forcing and its design sensitivity propagate to the response, response statistics, and optimization gradients.
\end{itemize}
Because the surrogate represents the localized nonlinear forcing rather than an output-specific response map, a single trained model can support different response quantities, statistical measures, and design criteria associated with the same dynamical system. 

The method is demonstrated on two nonlinear dynamical design problems. The first considers the collector area design of a solar heating system with uncertain heat loss coefficients and nonlinear thermal dynamics \cite{kicsiny2014extended}. The second example addresses the design of the yield force of a hysteretic isolation layer in a 100 degree of freedom building subjected to earthquake excitation \cite{wen1976method,kamalzare2015efficient,de2017efficient}. Similar optimization problems arise in friction-based isolation systems, for which the preferred design parameters may depend on ground motion characteristics and frictional heating \cite{li2024optimum}. In both examples, the neural network is used only for the localized nonlinear contribution, while the remaining system dynamics are evaluated from the nominal state space model. 

The rest of the paper is organized as follows. The design optimization problem under uncertainty is formulated in the next section. 
Thereafter, the accuracy of the proposed approach is analyzed through two propositions and three corollaries that relate errors in the neural network approximation of the pseudoforce and its design derivatives to errors in the response, mean and variance quantities, and their design gradients. 
Next, results from the two numerical examples are analyzed. Finally, the paper is concluded with a discussion of future extensions and applications of the proposed method.

\section{Problem Formulation} \label{sec:prob_form} 
In this section, we define the optimization problem for design under uncertainty. 
The governing differential equation of the dynamical system that contains the nonlinear uncertain components can be given in the state space form as follows (following Gaurav et al.  \cite{wojtkiewicz2011efficient} and De et al.  \cite{de2018computationally}),
\begin{equation}\label{eq:state-space}
\begin{split}
\dot{\Xm}(t;\thetaa,\xii) &= \Am \Xm(t;\thetaa,\xii) + \Bm \wm(t) + \Lm \gm\left(\Xmb(t;\thetaa,\xii);\thetaa,\xii\right), \quad \Xm(0) = \xm_0;\\
\Ym(t;\thetaa,\xii) &= \Cm \Xm(t;\thetaa,\xii) + \Dm \wm(t), \\
\end{split}
\end{equation}
where $\Xm(t;\thetaa,\xii)\in\mathbb{R}^{n}$ is the state vector; $\Am\in \mathbb{R}^{n\times n}$ is the state matrix; $\wm(t)\in\mathbb{R}^{\nw}$ is the external force vector; $\Bm\in\mathbb{R}^{n\times \nw}$ is the influence matrix for $\wm(t)$;  $\gm(\cdot;\cdot,\cdot)\in\mathbb{R}^{\ngm}$ arises from the nonlinear and uncertain components in the system and is a nonlinear function of a subset of the state, \textit{i.e.}, $\Xmb(t;\thetaa,\xii)=\Gm\Xm(t;\thetaa,\xii)$ with $\Gm\in\mathbb{R}^{\no\times n}$, $\no\ll n$, the uncertain parameter $\xii$ with known probability distribution $p_{\!_{\bm\Xi}}(\xii)$, and design parameter $\thetaa\in\bm\Theta$; $\Lm\in\mathbb{R}^{n\times \ngm}$ is the influence matrix for the nonlinear function $\gm(\cdot;\cdot,\cdot)$; and $\xm_0$ is the initial state vector. The output of the nonlinear dynamical system is $\Ym(t;\thetaa,\xii)\in\mathbb{R}^{\ny}$ with influence matrices as $\Cm\in\mathbb{R}^{\ny\times n}$ and $\Dm\in\mathbb{R}^{\ny\times \nw}$ for the state vector $\Xm(t;\thetaa,\xii)$ and external force $\wm(t)$, respectively. Note that the dynamical system can also contain nonlinearity and uncertainty in components other than the designed components. This can be easily handled in the above state space formulation, where all the terms that contain nonlinearity and uncertainty are collected in $\gm(\Xmb(t;\thetaa,\xii);\thetaa,\xii)$. 

A robust design of components of the nonlinear dynamical system is performed in this paper, where a combination of the expected value and variance of the normalized cost function $J(\Ym(t;\thetaa,\xii))/J_0$ and feasibility condition $h_j(\Ym(t;\thetaa,\xii))/h_{0j}$ for $j=1,\dots,n_h$ are used as the objective and constraints, respectively, where $J_0>0$ and $h_{0j}$ are fixed reference scales that are independent of $\thetaa$ and $\xii$. 
Hence, the optimization problem becomes 
\begin{equation} \label{eq:opt_def}
\begin{split}
    &\min_{\thetaa\in\bm\Theta}~~  \Exp_{\xii}\left[\frac{1}{J_0}J(\Ym(t;\thetaa,\xii))\right]+\lambda_J\mathrm{Var}\left(\frac{1}{J_0}J(\Ym(t;\thetaa,\xii))\right)\\
    &\text{subject to   } \Exp_{\xii}\left[\frac{1}{h_{0j}}h_j(\Ym(t;\thetaa,\xii))\right] + \lambda_{h_j} \mathrm{Var}\left(\frac{1}{h_{0j}}h_j(\Ym(t;\thetaa,\xii))\right) \leq 0;\quad j=1,\dots,n_h, 
\end{split}
\end{equation}
where $\mathrm{Var}_{\xii}(\cdot)$ denotes the variance of its argument with respect to the probability density $p_{\!_{\bm\Xi}}(\xii)$; $\lambda_J$ and $\lambda_{h_j}$ for $j=1,\dots,n_h$ are non-dimensional constants that show the importance of the variance in objective and constraints. 
Here variance is used as the dispersion measure following Beyer and Sendhoff \cite{beyer2007robust} and De et al.~\cite{de2020topology}, since it provides a smooth quadratic penalty. 
Using $n_s$ samples  $\{\xii_i\}_{i=1}^{n_s}$ for the uncertain parameters from the probability distribution $p_{\boldsymbol\Xi}(\xii)$, 
the above design optimization problem can be restated as 
\begin{equation} \label{eq:opt_def_ns} 
\begin{split}
    &\min_{\thetaa\in\bm\Theta}~~
    \frac{\overline{J}}{J_0} + \frac{\lambda_J}{J_0^2(n_s-1)} \sum_{i=1}^{n_s}\left[ J(\Ym(t;\thetaa,\xii_i)) - \overline{J}   \right]^2 \\
    &\text{subject to   }
    \frac{\overline{h}_j}{h_{0j}} + \frac{\lambda_{h_j}}{h_{0j}^2(n_s-1)} \sum_{i=1}^{n_s} \left[h_j(\Ym(t;\thetaa,\xii_{i})) - \overline{h}_j \right]^2\leq 0 ;\quad j=1,\dots,n_h, 
\end{split}
\end{equation}
where $\overline{J}:=\frac{1}{n_s}\sum_{i=1}^{n_s}J(\Ym(t;\thetaa,\xii_i))$ and $\overline{h}_j:=\frac{1}{n_s}\sum_{i=1}^{n_s}h_j(\Ym(t;\thetaa,\xii_i))$ for $j=1,\dots,n_{h}$ are corresponding average values estimated using $n_s$ samples. 
However, $n_s$ needs to be large for accurate estimation of the expectation or variance for most of the problems. In this paper, a sequential quadratic programming is used to solve the above optimization problem \cite{nocedal2006numerical}, which uses gradient information to construct a quadratic subproblem at every iteration.


\section{Proposed Methodology} 
The proposed use of neural networks for efficient calculation of the response and its sensitivity with respect to the design parameters is discussed herein. Thereafter, the accuracy of this proposed approach is discussed, where two propositions are presented that connect the error in the neural network modeling of the pseudoforce components to the objective, constraints, and their gradients. 

\begin{figure}[!htb]
    \centering
    \includegraphics[width=0.9\linewidth]{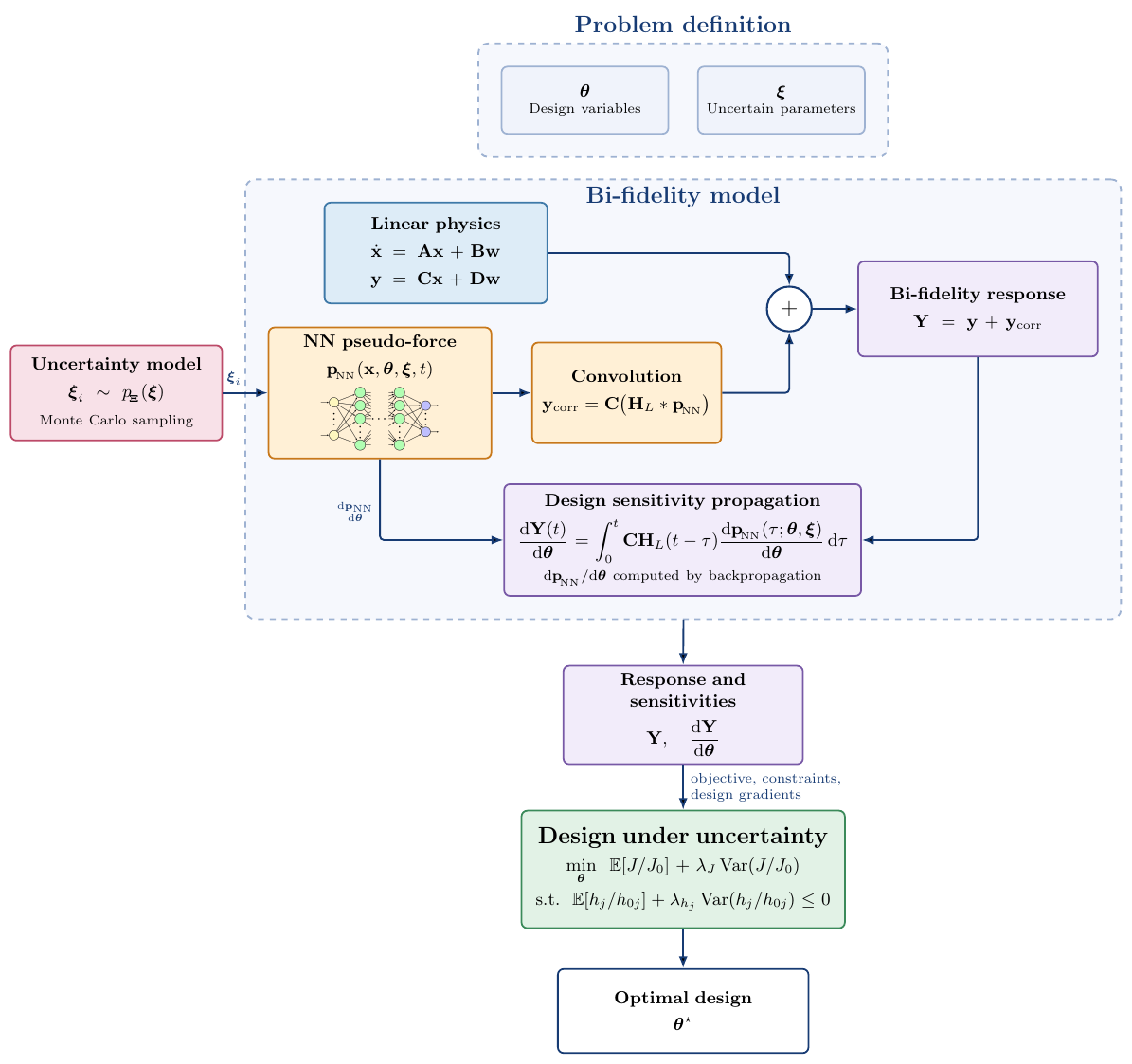}
    \caption{Steps of the proposed bi-fidelity approach of nonlinear dynamical system design under uncertainty.}
    \label{fig:schem}
\end{figure}

\subsection{Bi-fidelity Response and Sensitivity Calculations} \label{sec:eff_res}

The response of the uncertain nonlinear dynamical system in \eqref{eq:state-space} can be written as the summation of the {\textit{low-fidelity}} response $\xm(t)$ from a nominal linear system and a correction term $\xnl(t;\thetaa,\xii)$ that incorporates the effects of the nonlinearity, uncertainty, and the design component in the system as follows 
\begin{equation}
\Xm(t;\thetaa,\xii) = \xm(t)+\xnl(t;\thetaa,\xii), 
\end{equation} 
where the state $\xm(t)$ and the output $\ym(t)$ of the {\textit{low-fidelity}} nominal linear system corresponding to the nonlinear and uncertain system in \eqref{eq:state-space} are given by 
\begin{equation}\label{eq:det_response}
\begin{split}
\dot{\xm}(t) &= \Am \xm(t) + \Bm \wm(t), \quad \xm(0)=\xm_0;\\
\ym(t) &= \Cm \xm(t) + \Dm \wm(t).
\end{split}
\end{equation}
Hence, to estimate the response of the nominal linear system the impulse response function $\Hb(t) = \exp(\Am t) \Bm$ for the external force can be used as 
\begin{equation}\label{eq:linresp}
\begin{split}
\xm(t) &= \exp(\Am t) \xm_0 + \int_0^t \Hb(t-\tau)\wm(\tau)\drm\tau;\\
\ym(t) &= \Cm \exp(\Am t) \xm_0 + \int_0^t \Cm\Hb(t-\tau)\wm(\tau)\drm\tau + \Dm \wm(t).\\
\end{split}
\end{equation}
The estimation of $\xnl(t;\thetaa,\xii)$, however, provides challenges as it implicitly depends on the state vector $\Xmb(t;\thetaa,\xii)$. Herein, neural networks are used to efficiently estimate  $\xnl(t;\thetaa,\xii)$ by modeling a pseudoforce as
$\ppm(t;\thetaa,\xii) \approx \gm\left( \Xmb(t;\thetaa,\xii);\thetaa,\xii \right)$. 
For the neural networks, the design parameters $\thetaa$, the uncertain parameters $\xii$, and a subset of the nominal linear response $\xo(t) = \Gm \xm(t)$ are used as inputs. 
The correction term can then be given by 
\begin{equation}\label{eq:xnl}
\begin{split}
\xnl(t;\thetaa,\xii) &= \int_0^t \Hl (t-\tau) \ppm\left(\tau;\thetaa,\xii\right) \drm \tau,\\
\end{split}
\end{equation}
where the impulse response function $\Hl(t) = \exp(\Am t ) \Lm$ is due to the nonlinear pseudoforce. 
Once the correction term is estimated the output of \eqref{eq:state-space} becomes
\begin{equation}
\begin{split}
\Ym(t;\thetaa,\xii) &= \ym(t) +  
\Cm\xnl(t;\thetaa,\xii) \\
&= \ym(t) +  
\int_0^t \Cm\Hl (t-\tau) \ppm\left(\tau;\thetaa,\xii\right) \drm \tau. 
\end{split}
\end{equation} 
In the present formulation, the nominal matrices $\Am$, $\Cm$, and $\Lm$, the initial conditions, and the external input are independent of the design variables. Therefore, the design dependence enters through the pseudoforce. 
The sensitivity of the output $\Ym(t;\thetaa,\xii)$ is estimated as follows \cite{wojtkiewicz2014efficient} 
\begin{equation}
    \frac{\drm \Ym(t;\thetaa,\xii)}{\drm \thetaa} = \int_0^t \Cm \Hl(t-\tau) \frac{\drm \ppm\left(\tau;\thetaa,\xii\right)}{\drm \thetaa} \drm \tau
\end{equation}
Herein, the sensitivity of the pseudoforce is estimated using backpropagation through the neural network that is used to model the pseudoforce $\ppm(t;\thetaa,\xii)$.
 For a time discretization $t \in \{k\Delta t\}_{k=0}^{n_t}$, the output $\Ym(t;\thetaa,\xii)$ and the output sensitivity ${\drm \Ym(k\Delta t;\thetaa,\xii)}/{\drm \thetaa}$ with respect to the design parameters $\thetaa$ can be written as
\begin{equation} \label{eq:state_space_op}
    \begin{split}
        \Ym(k\Delta t;\thetaa,\xii) &= \ym(k\Delta t) + \Delta t \sum_{i=0}^{k-1} \Cm\Hl(k\Delta t - i\Delta t) \ppm(i\Delta t;\thetaa,\xii); \\
        \frac{\drm \Ym(k\Delta t;\thetaa,\xii)}{\drm \thetaa} &= \Delta t \sum_{i=0}^{k-1}  \Cm \Hl(k\Delta t-i\Delta t) \frac{\drm \ppm\left(i\Delta t;\thetaa,\xii\right)}{\drm \thetaa}; \quad k=0,1,\dots,n_t,\\ 
    \end{split} 
\end{equation} 
with $\Ym(0)=\ym(0)$. 
Note that the corrective term $\xnl(t;\thetaa,\xii)$ can be efficiently obtained from $\ppm(t;\thetaa,\xii)$ by performing the convolution using FFT. A schematic of the proposed approach is shown in Figure \ref{fig:schem}. 

\textbf{Remarks:} 
The following remarks can be made about the approach here. 
\\(1) The pseudoforce $\ppm(t;\thetaa,\xii)$ can also be estimated by solving a nonlinear Volterra integral equation of the second kind in a non-standard form given by \cite{wojtkiewicz2011efficient,de2018computationally} 
\begin{equation}\label{eq:nvie}
\begin{split}
&\ppm(t;\thetaa,\xii) - \gm\left(\xmb(t)+\int_0^t \Hl (t-\tau) \ppm\left(\tau;\thetaa,\xii\right) \drm \tau;\thetaa,\xii\right)=\bm{0}.\\
\end{split}
\end{equation}
The sensitivity of the pseudoforce with respect to the design parameter $\theta$ can be estimated as 
\begin{equation} \label{eq:nvie_sens}
    \frac{\drm \ppm(t;\thetaa,\xii)}{\drm \theta} = \frac{\partial \gm\left( \Xmb(t;\thetaa,\xii);\thetaa,\xii \right)}{\partial \theta} + \frac{\partial \gm\left( \Xmb(t;\thetaa,\xii);\thetaa,\xii \right)}{\partial \Xm(t;\thetaa,\xii)}\frac{\drm \Xm(t;\thetaa,\xii)}{\drm \theta}. 
\end{equation} 
\eqref{eq:nvie} and \eqref{eq:nvie_sens} can be solved to obtain the response and its sensitivity but it requires an iterative strategy due to the nonlinearity \cite{wojtkiewicz2011efficient,de2018computationally,de2017efficient}. 
\\
(2) The dimension $\ngm$ of  $\gm\left( \Xmb(t;\thetaa,\xii);\thetaa,\xii \right)$ can be high if many of the dynamical system's components contain nonlinearity or uncertainty. As a result, the number of outputs of the neural network required to train increases since each output corresponds to each dimension of $\gm\left( \Xmb(t;\thetaa,\xii);\thetaa,\xii \right)$. This may increase the offline cost of training the network; however,  during the optimization process, the cost of evaluating the neural networks for different realizations of the uncertain parameters $\xii$ remains small. 

\subsection{Datasets Used in Training of the Neural Networks} \label{sec:datasets} 

Two datasets are used to train and test the neural networks. To construct these datasets, i.i.d. samples of the uncertain parameters are generated from their respective probability distributions whereas samples from the design parameters are randomly selected from their domain assuming a uniform distribution between their lower and upper limits. 
For these realizations of the uncertain and design parameters, two datasets, namely, the training dataset $\Dtr=\Bigg\{\Big\{\ppm^{}\left( t_i;\thetaa_j^{(\mathrm{tr})},\xii_j^{(\mathrm{tr})}\right)\Big\}_{i=1}^{n_t},\thetaa_j^{(\mathrm{tr})},\xii_j^{(\mathrm{tr})}\Bigg\}_{j=1}^{\Ntr}$ and the validation dataset $\Dval =  \Bigg\{\Big\{\ppm^{}\left(t_i;\thetaa_j^{(\mathrm{val})},\xii_j^{(\mathrm{val})}\right)\Big\}_{i=1}^{n_t},\thetaa_j^{(\mathrm{val})},\xii_j^{\mathrm{(val)}}\Bigg\}_{j=1}^{\Nval}$ are generated by solving \eqref{eq:nvie}, where $\left\{\thetaa_j^{(\mathrm{tr})},\xii_j^{(\mathrm{tr})}\right\}_{j=1}^{\Ntr}$ and $\left\{\thetaa_j^{(\mathrm{val})},\xii_j^{(\mathrm{val})}\right\}_{j=1}^{\Nval}$ do not overlap. During training, the gradients are estimated using the training dataset $\Dtr$, the relative validation root mean squared error (RMSE) is defined using the validation dataset $\Dval$ as 
\begin{equation}\label{eq:val_rmse} 
\text{Relative validation RMSE} = \frac{\sum_{j=1}^{N_\mathrm{val}}\left[\sum_{k=0}^{n_t}\Big\lVert \ppm_{}^{}\left(t_k;\thetaa_j^{(\mathrm{val})},\xii_j^{(\mathrm{val})}\right) - \pnn^{}\left(t_k;\thetaa_j^{(\mathrm{val})},\xii_j^{(\mathrm{val})}\right) \Big\rVert_2^2\right]^{1/2}}{\sum_{j=1}^{N_\mathrm{val}}\left[\sum_{k=0}^{n_t}\Big\lVert \ppm\left(t_k;\thetaa_j^{(\mathrm{val})},\xii_j^{(\mathrm{val})}\right) \Big \rVert_2^2\right]^{1/2}},
\end{equation} 
where the subscript ``NN'' is used for predictions from the neural network, $t_k=k\Delta t$, 
and $\lVert \cdot \rVert_2$ is the Euclidean norm of its argument. 
We select the network complexity following De et al.~\cite{de2020transfer}. 
This training procedure takes a few hours on a modern desktop. However, once trained these networks can be used to efficiently calculate the response of the dynamical system for different realizations of the design and uncertain parameters. Further, the same network can be used for optimizing different objectives as the function $\gm\left(\Xmb(t;\thetaa,\xii);\thetaa,\xii\right)$ is modeled using the neural network instead of the cost function $J(\Ym(t;\thetaa,\xii))$. Next, the accuracy of the proposed approach is established.

\subsection{Accuracy of the Bi-fidelity Approximation} 

In this section, upper bounds are derived for the approximation errors in the objective, constraints, and their gradients under suitable regularity and moment conditions. These bounds are expressed in terms of the neural network approximation errors in the pseudoforce and its design derivatives, thereby identifying the quantities that must be controlled to improve the accuracy of the proposed bi-fidelity model. However, time discretization errors are not considered, since the analysis focuses exclusively on the error introduced by approximating the pseudoforce $\ppm(t;\thetaa,\xii)$ with a neural network. 

The analysis is stated for a generic scalar response functional $R(\mathcal{Y}(\thetaa,\xii))$, which may represent either the objective function or a constraint. The matrix $\mathcal{Y}(\thetaa,\xii)$ contains the selected system outputs $\Ym(t;\thetaa,\xii)$ evaluated for a realization of the design and uncertain parameters at the time instances $\mathcal{T}\subseteq\{k\Delta t\}_{k=0}^{n_t}$, where $t\in[0,t_{\max}]$, $t_{\max}=n_t\Delta t$, and the time domain is discretized as $\{k\Delta t\}_{k=0}^{n_t}$. In the following, the matrix $\mathcal{Y}(\thetaa,\xii)$ is vectorized as $\vecy$, \textit{i.e.}, the vector $\vecy$ denotes the output histories stacked over the time instances in $\mathcal T$. 
A subscript ``NN'' denotes a quantity obtained using the neural network approximation of the  pseudoforce. In the rest of this section, the explicit dependence of each quantity on $\thetaa$ and $\xii$ are omitted for brevity.

For a scalar-, vector-, or matrix-valued random quantity $q(\xii)$, define the norm $\|q\|_{L^r(p_{\boldsymbol{\Xi}})}:=\left(\Exp_{\xii}\left[\|q(\xii)\|^r\right]\right)^{1/r},~ r\geq 1$, where the norm inside the expectation is the absolute value for a scalar, the Euclidean norm for a vector, and the Frobenius norm for a matrix. 
Throughout this analysis, the nominal states, the output matrix $\Cm$, and the impulse response operator $\Hm_{L,k}$ are independent of the design variables. Hence, the design dependence of the output enters only through the pseudoforce, and 
\begin{equation}
\frac{\drm\Ym_k}{\drm\thetaa} = \Delta t \sum_{i=0}^{k} \Cm\Hm_{L,k-i} \frac{\drm\ppm_i}{\drm\thetaa}. 
\end{equation} 
At $t_i=i\Delta t$, define the pseudoforce error as $\varepsilon_i=\left\|\ppm_i-\ppm_{i,\mathrm{NN}}\right\|_2$, and the pseudoforce-sensitivity error as $\varepsilon_i^S=\left\|\frac{\drm \ppm_i}{\drm\thetaa}-\frac{\drm \ppm_{i,\mathrm{NN}}}{\drm\thetaa}\right\|_F$.

\begin{theorem} \label{prop1} 
Suppose that ${R}(\cdot)$ is Lipschitz continuous with respect to the stacked output history, \textit{i.e.},
\begin{equation}
    \left|R(\mathcal{Y})-R(\mathcal{Y}')\right|\leq\Lambda_R\left\|\operatorname{vec}(\mathcal{Y})-\operatorname{vec}(\mathcal{Y}')\right\|_2
\end{equation}
for all admissible output histories $\mathcal{Y}$ and $\mathcal{Y}'$. Assume further that $\|\Cm\|_2\leq\Gamma_C,~\|\Hm_{L,k}\|_2\leq\Gamma_H,~k=0,\ldots,n_t$. Then
\begin{equation}\label{eq:prop1}
    \begin{split}
        &\left\|R(\mathcal{Y})-R(\mathcal{Y}_{\mathrm{NN}})\right\|_{L^2(p_{\boldsymbol{\Xi}})}\leq\Delta t\Lambda_R\Gamma_C\Gamma_H\sum_{t_k\in\mathcal T}\sum_{i=0}^{k}\left\|\varepsilon_i\right\|_{L^2(p_{\boldsymbol{\Xi}})}.
    \end{split}
\end{equation} 


\end{theorem} 

A proof of this proposition is given in Appendix \ref{sec:proof_p1}. 

\textbf{Remarks:}\\
(1) As $\Lambda_R$, $\Gamma_C$, and $\Gamma_H$ depend on the expression of objective or constraint and the state space formulation of the problem, the error due to the use of a neural network only enters through $\epsilon_i$ in the bound \eqref{eq:prop1}. This relation also shows that the error could accumulate in the latter part of the time history. 
Here the first summation is over $t\in\mathcal{T}$. Hence, if the quantity $R(\cdot)$ only depends on the response for one time instance the need for this summation is removed. 
The second summation is due to the convolution of the impulse response $\Hl(t)$ and the nonlinear pseudoforce $\ppm(t;\thetaa,\xii)$.   
To avoid this error for long time histories, a better network (possibly larger) or a family of smaller but accurate networks to model different parts of the time history need to be trained. Note that there are many improvements proposed in the literature for training networks for long time-histories, such as recurrent neural networks \cite{goodfellow2016deep} or networks with attention mechanism \cite{vaswani2017attention}, and so on. However, a detailed comparison of these different architectures of the neural network is beyond the scope of the current paper. 
In the numerical examples studied in this paper, a reasonably large neural network is shown to provide accurate predictions for $\ppm(t;\thetaa,\xii)$. 

\noindent(2) During implementation, the neural network should be validated on the pseudoforce rather than only on the final objective or constraint values. Since the bounds depend directly on \(\epsilon_i=\|\mathbf{p}(i\Delta t)-\mathbf{p}_{\mathrm{NN}}(i\Delta t)\|_2\), the training and validation workflow should track this error over time, over uncertainty samples, and over the design domain. The surrogate should be tested more carefully at later portions of the response history and at time instants that strongly influence the objective or constraints. 

\textbf{Corollary 1.} 
Let $Z=R(\mathcal{Y})/R_0,~Z_{\mathrm{NN}}=R(\mathcal{Y}_{\mathrm{NN}})/R_0$, where $R_0>0$ is a reference scale,
and define 
\begin{equation}
    \begin{split}
        \mathcal R(\thetaa) & = \Exp_{\xii}[Z] + \lambda_R \operatorname{Var}_{\xii}(Z), \\
        \mathcal R_{\mathrm{NN}}(\thetaa) & = \Exp_{\xii}[Z_{\mathrm{NN}}] + \lambda_R\operatorname{Var}_{\xii}(Z_{\mathrm{NN}})\\ 
    \end{split}
\end{equation} 
If \(Z,Z_{\mathrm{NN}}\in L^2(p_{\boldsymbol{\Xi}})\), then
\begin{equation}
    \left|\mathcal R(\thetaa)-\mathcal R_{\mathrm{NN}}(\thetaa)\right|\leq\left[1+\lambda_R\left(\sigma_Z+\sigma_{Z_{\mathrm{NN}}}\right)\right]\delta_R(\thetaa),
\end{equation}

where $\sigma_Z=\sqrt{\operatorname{Var}_{\xii}(Z)}$, $\sigma_{Z_{\mathrm{NN}}}=\sqrt{\operatorname{Var}_{\xii}(Z_{\mathrm{NN}})}$, and $\delta_R(\boldsymbol{\theta}):=\frac{1}{R_0}\Delta t \Lambda_R\Gamma_C\Gamma_H\sum_{t_k\in\mathcal T}\sum_{i=0}^{k}\left\|\varepsilon_i\right\|_{L^2(p_{\boldsymbol{\Xi}})}$. A proof of this corollary is given in Appendix \ref{sec:proof_c1}.

\begin{theorem} \label{proop2} 

Suppose that the gradient of $R$ with respect to the stacked output history is Lipschitz continuous, \textit{i.e.},
\begin{equation}
    \left\|\nabla_{\vey}R(\mathcal{Y}) - \nabla_{\vey}R(\mathcal{Y}')\right\|_2 \leq \Lambda_T \left\| \vey - \mathrm{vec}(\mathcal{Y}')
\right\|_2
\end{equation}
for all admissible output histories $\mathcal{Y}$ and $\mathcal{Y}'$. Define the exact and approximate output sensitivities as 
$\Sm_Y = \frac{\drm\vey}{\drm\thetaa}$, $\Sm_{Y,\mathrm{NN}} = \frac{\drm\veyn}{\drm\thetaa}$.
Assume that, for each $\thetaa\in\Theta$, $\nabla_{\vey}R(\mathcal{Y})\in L^2(p_{\boldsymbol{\Xi}})$ and $\Sm_{Y,\mathrm{NN}}\in L^2(p_{\boldsymbol{\Xi}})$.
Let $B_1(\thetaa) = \left\| \nabla_{\vey} R(\mathcal{Y}) \right\|_{L^2(p_{\boldsymbol{\Xi}})}$ and $B_2(\thetaa) = \left\|\Sm_{Y,\mathrm{NN}}\right\|_{L^2(p_{\boldsymbol{\Xi}})}$.
Then 
\begin{equation} \label{eq:prop2_1}
\begin{split}
&\left\|\nabla_{\thetaa}R(\mathcal{Y}) - \nabla_{\thetaa}R(\mathcal{Y}_{\mathrm{NN}})\right\|_{L^1(p_{\boldsymbol{\Xi}})} \leq B_1(\thetaa) \Delta t \Gamma_C\Gamma_H \sum_{t_k\in\mathcal T} \sum_{i=0}^{k} \left\|\varepsilon_i^S\right\|_{L^2(p_{\boldsymbol{\Xi}})} + \Lambda_T B_2(\thetaa)\Delta t \Gamma_C\Gamma_H \sum_{t_k\in\mathcal T} \sum_{i=0}^{k} \left\|\varepsilon_i\right\|_{L^2(p_{\boldsymbol{\Xi}})}.
\end{split}
\end{equation}

\end{theorem}  
A proof of this proposition is given in Appendix \ref{sec:proof_p2}. 

\textbf{Remarks:}\\ 
(1) Proposition~2 identifies two sources of error in the estimated design gradient. The first term arises from the error in the derivative of the learned pseudoforce and therefore measures how accurately the neural network reproduces the sensitivity of the nonlinear correction with respect to the design variables. The second term arises from the pseudoforce value error, which perturbs the system output and thereby changes the gradient of the response functional with respect to that output. Accurate response prediction alone is therefore insufficient to guarantee accurate design gradients; the derivatives of the learned pseudoforce must also be accurate. The constants $\Gamma_C$ and $\Gamma_H$ describe the amplification of these local approximation errors by the output operator and the nominal linear dynamics, while $B_1$, $B_2$, and $\Lambda_T$ characterize the sensitivity of the selected response functional and its associated output sensitivities. 

\noindent(2) Proposition~2 bounds the expected magnitude of the gradient error, \textit{i.e.}, its $L^1(p_{\boldsymbol{\Xi}})$ norm, using finite second moment assumptions. A corresponding $L^2(p_{\boldsymbol{\Xi}})$ bound can also be established. However, because the gradient error contains products of random output sensitivities and response gradients, such a result requires stronger integrability conditions. In particular, finite fourth moments are  sufficient through H\"older's inequality, $\|fg\|_{L^2}\leq\|f\|_{L^4}\|g\|_{L^4}$.
The $L^1$ formulation is adopted here because it controls the mean gradient error under the weaker assumption of finite second moments. 

(3) Response functionals involving temporal extrema require a locally unique active time for Proposition~2. For example, if $R(\mathcal{Y})=\max_{t_k\in\mathcal{T}}f(\Ym_k)$ has a locally unique maximum at $t_{k^\ast}$, then $R(\mathcal{Y})=f(\Ym_{k^\ast})$ in a neighborhood and the gradient bound applies. If the extreme value occurs at multiple time instances, $R$ remains Lipschitz, so Proposition~1 still applies, but it may not be differentiable. Note that a smooth approximation may be used when global differentiability is required. 

(4) Proposition 2 also shows that accurate response prediction alone is insufficient for gradient-based optimization. Both the pseudoforce error $\epsilon_i$ and its design derivative error $\epsilon_i^S$ contribute to the gradient of the response function. Accordingly, the numerical examples validate both $\ppm$ and $\frac{\mathrm{d}\ppm}{\mathrm{d}\thetaa}$. 

\textbf{Assumption 1 (Derivative consistency).} 
Assume that the trained pseudoforce surrogate is consistent with respect to both its values and its design derivatives, and that there exists a finite constant $\Lambda_p$ such that $\sum_{t_k\in\mathcal T}\sum_{i=0}^{k}\left\|\varepsilon_i^S\right\|_{L^2(p_{\boldsymbol{\Xi}})}\leq\Lambda_p\sum_{t_k\in\mathcal T}\sum_{i=0}^{k}\left\|\varepsilon_i\right\|_{L^2(p_{\boldsymbol{\Xi}})}$ for each $\thetaa\in\Theta$.

\textbf{Corollary 2.}
Under the conditions of Proposition~2 and Assumption~1,

\begin{equation}
\left\|\nabla_{\thetaa}R(\mathcal{Y})-\nabla_{\thetaa}R(\mathcal{Y}_{\mathrm{NN}})\right\|_{L^1(p_{\boldsymbol{\Xi}})}
\leq\left[\Lambda_pB_1(\thetaa)+ \Lambda_T B_2(\thetaa)\right]\Delta t \Gamma_C\Gamma_H\sum_{t_k\in\mathcal T}\sum_{i=0}^{k}\left\|\varepsilon_i\right\|_{L^2(p_{\boldsymbol{\Xi}})}.
\end{equation} 
A proof of this corollary is given in Appendix \ref{sec:proof_c2}.

\textbf{Remark:} 
Assumption~1 expresses first-order derivative consistency of the pseudoforce surrogate over the design and uncertainty domains. It is motivated by approximation in a weighted Sobolev space associated with $p_{\boldsymbol{\Xi}}$, in which the accuracy of both the pseudoforce and its design derivatives is considered. In the numerical examples, this consistency is assessed on an independent validation dataset by comparing the design derivatives obtained through automatic differentiation with reference sensitivities evaluated using central differences. 


\textbf{Corollary 3.}
Define $\gm = \nabla_{\thetaa}Z$, $\gm_{\mathrm{NN}} = \nabla_{\thetaa}Z_{\mathrm{NN}}$. 
Suppose differentiation and expectation can be interchanged and $Z,Z_{\mathrm{NN}}\in L^2(p_{\boldsymbol{\Xi}})$, $\gm,\gm_{\mathrm{NN}}\in L^2(p_{\boldsymbol{\Xi}})$, $\nabla_{\vey}R(\mathcal{Y})\in L^4(p_{\boldsymbol{\Xi}})$, $\Sm_{Y,\mathrm{NN}}\in L^4(p_{\boldsymbol{\Xi}})$, and $\varepsilon_i,\varepsilon_i^S\in L^4(p_{\boldsymbol{\Xi}}), ~ i=0,\ldots,n_t$.  
Then the $L^2(p_{\boldsymbol{\Xi}})$ form of the argument used in Proposition~2, along with H\"older's inequality, gives 
\begin{equation}
    \left\|\gm-\gm_{\mathrm{NN}}\right\|_{L^2(p_{\boldsymbol{\Xi}})}\leq\delta_R^S(\thetaa), 
\end{equation}
where 
\begin{equation}
    \delta_R^S(\thetaa)=\frac{1}{R_0}B_1^{(4)}(\thetaa) \Delta t \Gamma_C\Gamma_H \sum_{t_k\in\mathcal T} \sum_{i=0}^{k} \left\|\varepsilon_i^S\right\|_{L^4(p_{\boldsymbol{\Xi}})}+ \frac{1}{R_0}\Lambda_T B_2^{(4)}(\thetaa) \Delta t \Gamma_C\Gamma_H \sum_{t_k\in\mathcal T} \sum_{i=0}^{k} \left\|\varepsilon_i\right\|_{L^4(p_{\boldsymbol{\Xi}})}
\end{equation} 
with $B_1^{(4)}(\thetaa) = \left\|\nabla_{\vey} R(\mathcal{Y})\right\|_{L^4(p_{\boldsymbol{\Xi}})}$ and $B_2^{(4)}(\thetaa) = \left\|\Sm_{Y,\mathrm{NN}}\right\|_{L^4(p_{\boldsymbol{\Xi}})}$. 
Then 
\begin{equation}
    \left\|\nabla_{\thetaa}\mathcal R(\thetaa) - \nabla_{\thetaa}\mathcal R_{\mathrm{NN}}(\thetaa)\right\|_2 \leq \left(1+2\lambda_R\sigma_{Z_{\mathrm{NN}}}\right)\delta_R^S(\thetaa) + 2\lambda_R\sigma_{\gm}\delta_R(\thetaa)
\end{equation}
where $\sigma_{\gm} = \left\|\gm - \Exp_{\xii}[\gm]\right\|_{L^2(p_{\boldsymbol{\Xi}})}$. A proof of this corollary is given in Appendix \ref{sec:proof_c3}. 


\textbf{Remark:}
Corollary~3 provides a bound for the gradient of the complete mean-variance objective or constraint. Unlike Proposition~2, which bounds the expected magnitude of the gradient error using only second moments, this result requires an $L^2$ bound on the gradient error.
Because the gradient error contains products of random output sensitivities and response gradients, sufficient conditions for this stronger result include finite fourth moments of these quantities and of the pseudoforce value and sensitivity errors.

\section{Numerical Examples} 
Two numerical examples utilizing nonlinear dynamical systems are used in this section to illustrate the proposed approach. PyTorch \cite{paszke2017automatic} is used to implement the neural networks for the examples. The required sensitivities are computed through PyTorch automatic differentiation \cite{paszke2017automatic} and propagated through the linear state space model to obtain the gradients of the state variables.
The constrained optimization problem is solved using the Sequential Least Squares Programming (SLSQP) algorithm \cite{nocedal2006numerical} available in the SciPy optimization library \cite{virtanen2020scipy}.
At each SLSQP iteration, a new set of \(n_s\) independent and identically distributed (i.i.d.) uncertainty realizations is generated. The same sample set is used to evaluate the objective, constraints, and their gradients within that iteration. This resampling avoids optimizing repeatedly against a single finite-sample realization to avoid finite-sample bias, while the large batch size limits the resulting sampling noise.
A convergence tolerance (\texttt{ftol}) of $10^-6$ is used for the SLSQP optimizer. As described in the SciPy documentation \cite{scipy_slsqp_docs}, this tolerance governs convergence checks involving changes in the objective function, the optimization step size, the gradient of the Lagrangian, and the aggregate magnitude of constraint violations.
For brevity, the dependence of the quantities on $t$ and $\xii$ are omitted in places in this section. Herein, the accuracy of the estimates is measured using the relative root mean squared error (RMSE) given by $\frac{\lVert \widehat{\ym} - \ym \rVert_2}{\lVert \ym \rVert_2}$, where $\widehat{\ym}$ is the estimated quantity and $\ym$ is the true response. 

\subsection{Example I: Solar Heating System}

\begin{figure}[!htb]
    \centering
    \includegraphics[width=0.75\linewidth]{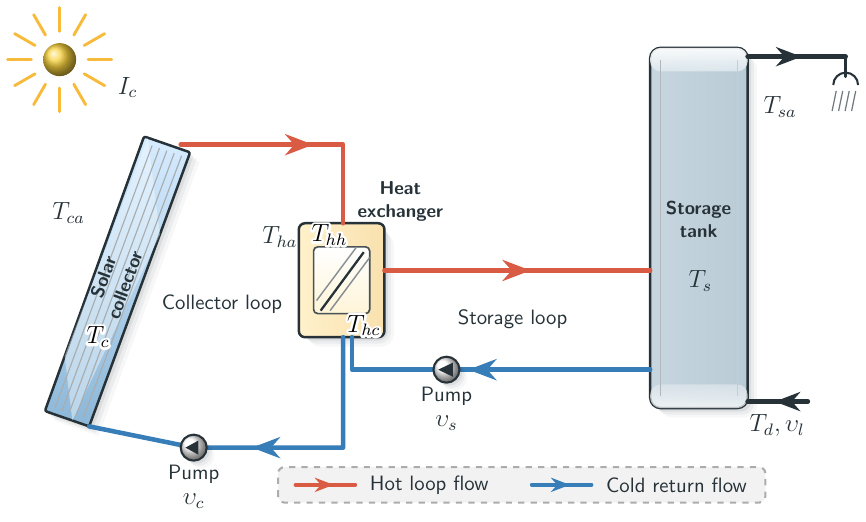}
    \caption{Schematic of the solar heating system used in Example I.}
    \label{fig:solar_schem}
\end{figure} 

\subsubsection{Description}
In the first example, we use a solar heating system  shown in Figure \ref{fig:solar_schem} and expressed through a system of four coupled ordinary differential equations of the first order \cite{kicsiny2014extended,buzas1998modelling,zeghib2011simulation}. The corresponding system of ordinary differential equations is 
\begin{subequations} \label{eq:solarODEs}
\begin{align}
    \dot{T}_c &= \frac{A_c\eta_0}{\rho_c c_c V_c}I_c + \frac{k_1 A_c}{\rho_c c_c V_c}(T_{ca} - T_c) - \frac{k_2 A_c}{\rho_c c_c V_c}(T_{ca} - T_c)^2 + \frac{\upsilon_c}{V_c}(T_{hh} - T_c); \\
    \dot{T}_{hh} &= \frac{\rho_c c_c \upsilon_c}{\frac{c_h m_h}{2}+\rho_c c_c \frac{V_h}{2}}(T_c - T_{hh}) + \frac{\varepsilon k_h A_h}{\frac{c_h m_h}{2}+\rho_c c_c \frac{V_h}{2}}(T_{hc} - T_{hh}) + \frac{\frac{A_a k_a}{2}}{\frac{c_h m_h}{2}+\rho_c c_c \frac{V_h}{2}}(T_{ha} - T_{hh});\\
    \dot{T}_{hc} &= \frac{\rho_s c_s \upsilon_s}{\frac{c_h m_h}{2}+\rho_s c_s \frac{V_h}{2}}(T_s - T_{hc}) + \frac{\varepsilon k_h A_h}{\frac{c_h m_h}{2}+\rho_s c_s \frac{V_h}{2}}(T_{hh} - T_{hc}) + \frac{\frac{A_a k_a}{2}}{\frac{c_h m_h}{2}+\rho_s c_s \frac{V_h}{2}}(T_{ha} - T_{hc});\\
    \dot{T}_s &= \frac{\upsilon_l}{V_s}(T_{d} - T_s) + \frac{\upsilon_s}{V_s}(T_{hc} - T_s) + \frac{A_s k_s}{\rho_s c_s V_s}(T_{sa} - T_s), 
\end{align}
\end{subequations}
where $T_c$ is the temperature of the solar collector, $T_{hh}$ is the temperature at the hot side of the heat exchanger (collector), $T_{hc}$ is the temperature at the cold side of the heat exchanger (storage), and $T_s$ is the temperature of the solar storage.  All temperatures have units of degrees Celsius.  The other time-dependent inputs are defined in Table~\ref{tab:var_time}, whereas  the constant parameters are listed in Table~\ref{tab:const_params} \cite{kicsiny2014extended}. The values of $\upsilon_c$ and $\upsilon_s$ are assumed as 0.98 m$^3$/h and 0.63 m$^3$/h, respectively, following Kicsiny et al.~\cite{kicsiny2014extended}, and the consumption load of the solar storage as $\upsilon_l=2\times10^{-4}$ m$^3$/s. 

\begin{table}[h!]
    \centering
    \caption{Time-dependent input variables}
    \label{tab:var_time}
    \begin{tabular}{ccc}
    \toprule
    Variable & Description & Units \\
    \midrule
    $I_c$ & global solar irradiance on collector surface & W/m$^2$ \\
    $T_{ca}$ & ambient temperature of the solar collector & \degree C \\
    $T_{ha}$ & ambient temperature of the heat exchanger & \degree C \\
    $T_{sa}$ &  ambient temperature of the solar storage & \degree C \\
    $T_d$ & {tap water temperature} & \degree C \\
    $\upsilon_c$ & {volumetric flow in the collector loop} & m$^3$/s \\
    $\upsilon_s$ & {volumetric flow in the storage loop} & m$^3$/s \\
    $\upsilon_l$ & {consumption load of the solar storage} & m$^3$/s\\
    \bottomrule
    \end{tabular}
\end{table}

\begin{table}[h!]
    \centering
        \caption{Constant model parameters}
    \label{tab:const_params}
    \begin{tabular}{cccc}
    \toprule
    Variable & Description & Value & Units \\
    \midrule
    $\eta_0$ & optical efficiency of the collector & 0.74 & -- \\
    $\rho_c$ & collector fluid density & 1034 & kg/m$^3$ \\
    $c_c$ & specific heat capacity of the collector fluid & 3623 & J/(kg K) \\
    $V_c$ & volume of the solar collector field & 0.027 & m$^3$ \\
    $c_h$ & specific heat capacity of the heat exchanger material & 464.76 & J/(kg K) \\
    $m_h$ & mass of the empty heat exchanger & 37 & kg \\
    $V_h$ & volume of the heat exchanger & 0.005 & m$^3$ \\
    $\varepsilon k_h$ & heat transfer coefficient inside the heat exchanger & 2461.5 & W/(m$^2$ K) \\
    $A_h$ & surface area of the heat transfer inside the heat exchanger & 2 & m$^2$ \\
    $k_a$ & heat loss coefficient of the heat exchanger to the ambiance & 0 & W/(m$^2$ K) \\
    $A_a$ & surface area of the heat exchanger to the ambiance & 0.24 & m$^2$ \\
    $\rho_s$ & density of water (storage fluid) & 1000 & kg/m$^3$ \\
    $c_s$ & specific heat capacity of water & 4200 & J/(kg K) \\
    $V_s$ & volume of the solar storage & 2 & m$^3$ \\
    $A_s$ & outside surface area of the solar storage & 4 & m$^2$ \\
    $k_s$ & heat loss coefficient of the solar storage to the ambiance & 1 & W/(m$^2$ K) \\
    \bottomrule
    \end{tabular}
\end{table}

The system of equations in \eqref{eq:solarODEs} can be reorganized into state space format according to the following state space vector
\begin{equation}
    \mathbf{x} = [T_c, T_{hh}, T_{hc}, T_s]\transpose
\end{equation}
which then allows for matrix $\Am$ to be defined as 
\begin{equation}
    \Am = \left[\begin{array}{cccc}
        - \frac{\upsilon_c}{V_c}& \frac{\upsilon_c}{V_c} & 0 & 0 \\
        \frac{\rho_c c_c \upsilon_c}{\frac{c_h m_h}{2}+\rho_c c_c \frac{V_h}{2}} & -\frac{\rho_c c_c \upsilon_c+\varepsilon k_h A_h+\frac{A_a k_a}{2}}{\frac{c_h m_h}{2}+\rho_c c_c \frac{V_h}{2}} 
        & \frac{\varepsilon k_h A_h}{\frac{c_h m_h}{2}+\rho_c c_c \frac{V_h}{2}} & 0 
        \\
        0 & \frac{\varepsilon k_h A_h}{\frac{c_h m_h}{2}+\rho_s c_s \frac{V_h}{2}} & -\frac{\rho_s c_s \upsilon_s+\varepsilon k_h A_h+\frac{A_a k_a}{2}}{\frac{c_h m_h}{2}+\rho_s c_s \frac{V_h}{2}} 
        & \frac{\rho_s c_s \upsilon_s}{\frac{c_h m_h}{2}+\rho_s c_s \frac{V_h}{2}} \\
        0 & 0 & \frac{\upsilon_s}{V_s} & -\frac{\upsilon_l}{V_s} - \frac{\upsilon_s}{V_s} - \frac{A_s k_s}{\rho_s c_s V_s}
    \end{array}\right]
\end{equation}

We assume the collector field surface area to be the design parameter, \textit{i.e.}, $\theta=A_c$. 
Further, the external inputs $\wm$ are given as 
    $\wm = [T_{ha}, T_{sa},T_d]\transpose
    $. 
    The temperatures $T_{ha}$ and $T_{sa}$ follow a prescribed sinusoidal profile, increasing from $30^\circ\mathrm{C}$ to $40^\circ\mathrm{C}$ over the 5 hours period as below:
\begin{equation}
    T_{ha}(t)=T_{sa}(t)
    =30+10\sin\left(\frac{2\pi t}{4t_f}\right),
    \qquad 0\leq t\leq t_f,
\end{equation}
with $t_f=5$ hours. Tap water temperature is $T_d = 23 ^\circ$C$ $. 
Then $\Bm$ is given as
\begin{equation}
    \Bm = \left[\begin{array}{ccccc}
     0 & 0 & 0\\
     \frac{\frac{A_a k_a}{2}}{\frac{c_h m_h}{2}+\rho_c c_c \frac{V_h}{2}}  & 0 & 0\\
     \frac{\frac{A_a k_a}{2}}{\frac{c_h m_h}{2}+\rho_s c_s \frac{V_h}{2}}  & 0 & 0\\
     0 & \frac{A_s k_s}{\rho_s c_s V_s} & \frac{\upsilon_l}{V_s}
    \end{array}\right]. 
\end{equation}
To account for uncertainty in the thermal performance of the solar collector, the first-order heat loss coefficient $k_1$, the second-order heat loss coefficient $k_2$ are modeled as uncertain parameters in the governing equations, \textit{i.e.}, $\xii = [k_1,k_2]\transpose$.
The thermal loss coefficients are modeled using lognormal probability distributions as shown in Table \ref{tab:uncertainty}. The nonlinear and uncertain pseudoforce is given by the term 
\begin{equation} \label{eq:solar_g}
    p(t;\theta,\xii)=
        \frac{k_1 A_c}{\rho_c c_c V_c}(T_{ca} - T_c) -\frac{k_2 A_c}{\rho_c c_c V_c}(T_{ca} - T_c)^2 + \frac{A_c\eta_0}{\rho_c c_c V_c} I_c 
\end{equation}
while influence matrix $\Lm$ is simply
$    \Lm = 
    \left[\begin{array}{cccc}
        1 &
        0 &
        0 &
        0 
    \end{array}\right]\transpose$.
The initial temperatures are specified as $T_c(0)=46.5^\circ\mathrm{C}$, $T_{hh}(0)=40.1^\circ\mathrm{C}$,
$T_{hc}(0)=38.85^\circ\mathrm{C}$, and $T_s(0)=30^\circ\mathrm{C}$, respectively.
\begin{table}[!htb]
\centering
\caption{Probability distributions of the uncertain parameters.}
\label{tab:uncertainty}
\begin{tabular}{lcccc}
\hline \Tstrut 
Parameter & Distribution & Mean & Standard deviation \Bstrut\\
\hline \Tstrut 
$k_1$ (W/m$^2$K)
& Lognormal
& 7.0
& 3.5 \\

$k_2$ (W/m$^2$K$^2$)
& Lognormal
& 0.027
& 0.0135 \Bstrut \\

\hline
\end{tabular}
\end{table}

\subsubsection{Optimization Problem}
The design optimization problem for this example uses an objective of maximizing the total thermal energy delivered from the heat exchanger to the storage tank using 
\begin{equation}
{J}(\theta,\xii) = \int_{0}^{t_{\max}} \rho_s c_s \upsilon_s\left(T_{hc}(t;\theta,\xii)-T_s(t;\theta,\boldsymbol{\xi})\right)\mathrm{d}t, 
\end{equation} 
where $t_{\max}$ denotes the maximum operating hours. 
To improve the numerical scaling of the optimization problem, the objective function is normalized by a reference value $J_0=10^8 $~J. 
In addition, three constraints are used while keeping the collector area $A_c$ as the design parameter $\theta$ within a bound $A_c\in [A_c^{\min}, A_c^{\max}]$. The first constraint specifies that the storage temperature $T_s$ must remain above a prescribed minimum delivery temperature $T_s^{\min}$ throughout the operating hours, 
\begin{equation}
{h}_1(\theta,{\xii})=T_s^{\min}-\min_{t\in[0,t_{\max}]}T_s(t;\theta,{\xii})\le0.
\end{equation}
This constraint is added so that sufficient thermal energy is always available for domestic hot water demand. 
To prevent overheating of the storage tank, the storage temperature is restricted by
\begin{equation}
{h}_2(\theta,{\xii})=\max_{t\in[0,t_{\max}]}T_s(t;\theta,{\xii})-T_s^{\max}\le0,
\end{equation}
where $T_s^{\max}$ denotes the maximum allowable storage temperature.
To avoid excessive collector temperatures, the collector temperature is constrained by
\begin{equation}
{h}_3(\theta,\boldsymbol{\xi})=\max_{t\in[0,t_{\max}]}T_c(t;\theta,\boldsymbol{\xi})-T_c^{\max}\le0,
\end{equation}
where $T_c^{\max}$ is the maximum allowable collector temperature. To provide consistent numerical scaling, each constraint is normalized by its corresponding prescribed temperature limit, \textit{i.e.}, $h_{01}=T_s^{\min}$, $h_{02}=T_s^{\max}$, and $h_{03}=T_c^{\max}$. 
Hence, the normalized robust optimization problem is given by 
\begin{equation}
    \begin{aligned}
    \min_{\theta}\quad
    & 
    \Exp_{\xii}
    \left[-\frac{1}{J_0}
    J(\theta,\xii)
    \right]
    +
    \lambda_J
    \operatorname{Var}
    \left[\frac{1}{J_0}
    J(\theta,\xii)
    \right]
    \\
    \text{subject to}\quad
    &
    \Exp_{\xii}
    \left[\frac{1}{h_{0j}}
    h_j(\theta,\xii)
    \right]
    +
    \lambda_{h_j}
    \operatorname{Var}
    \left[\frac{1}{h_{0j}}
    h_j(\theta,\xii)
    \right]
    \leq
    0,\quad j = 1,2,3 
    \\
    &
    A_c^{\min} \leq \theta \leq A_c^{\max},
    \end{aligned}
    \label{eq:op_ex1}
\end{equation} 
with parameter values used in the optimization problem listed in Table \ref{tab:opt_settings}.  



\subsubsection{Results} 
The pseudoforce is approximated using a feedforward neural network (FNN). The network consists of five inputs corresponding to the collector area $A_c$, the first- and second-order heat loss coefficients $k_1$ and $k_2$, the collector temperature obtained from the linear nominal system, and time $t$, followed by $N_H=4$ hidden layers with $m=100$ neurons each and a single output representing the pseudoforce. Hyperbolic tangent activation functions are employed. 
A total of $500$ Monte Carlo realizations are generated to construct the dataset, of which $450$ trajectories are used for training and the remaining $50$ trajectories are reserved for validation and performance evaluation, with a sampling interval of $\Delta t=30~\mathrm{s}$. 
The network parameters are optimized using the Adam optimizer \cite{kingma2014adam} with a learning rate of $4\times10^{-5}$ by minimizing the mean squared error (MSE) between the predicted and target pseudoforce values. The network is trained for 10,000 epochs. A step learning rate scheduler is employed, where the learning rate is reduced by a factor of $0.5$ every 1,000 epochs. The minimum relative validation error obtained during training is $2.464\times10^{-4}$ with a total training time of $3.33$ hours. 
\begin{figure}[!htb]
	\centering
	\begin{subfigure}[t]{0.45\textwidth}
		\centering
		\includegraphics[scale=0.45]{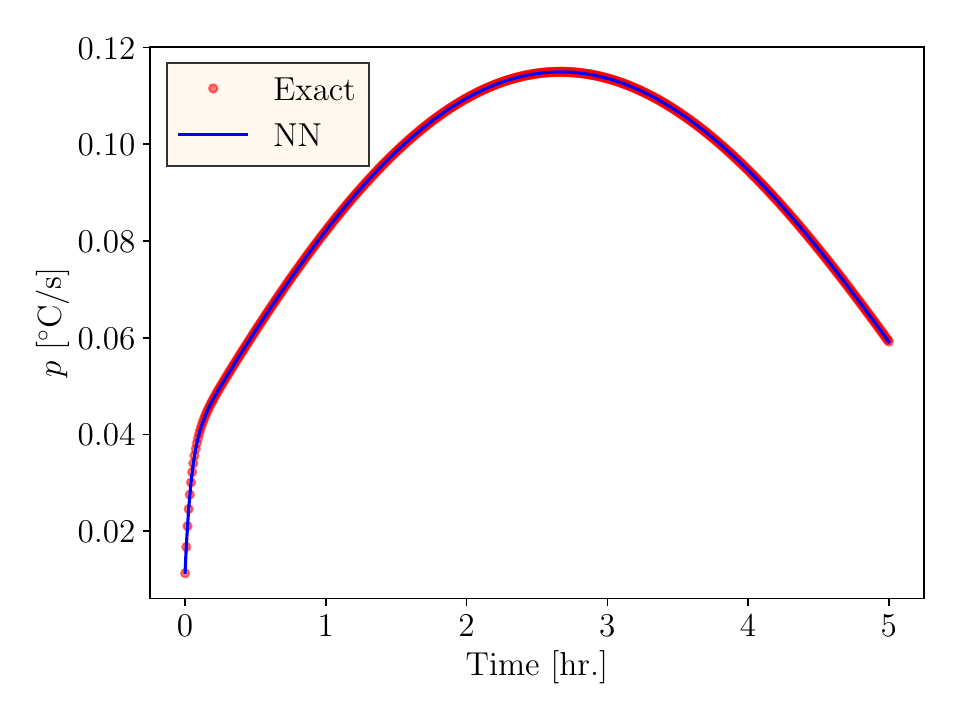}
		\caption{Comparison of estimated pseudoforce}\label{fig:solar_p_comp}
	\end{subfigure}
	\hfill
	\begin{subfigure}[t]{0.45\textwidth}
		\centering
		\includegraphics[scale=0.45]{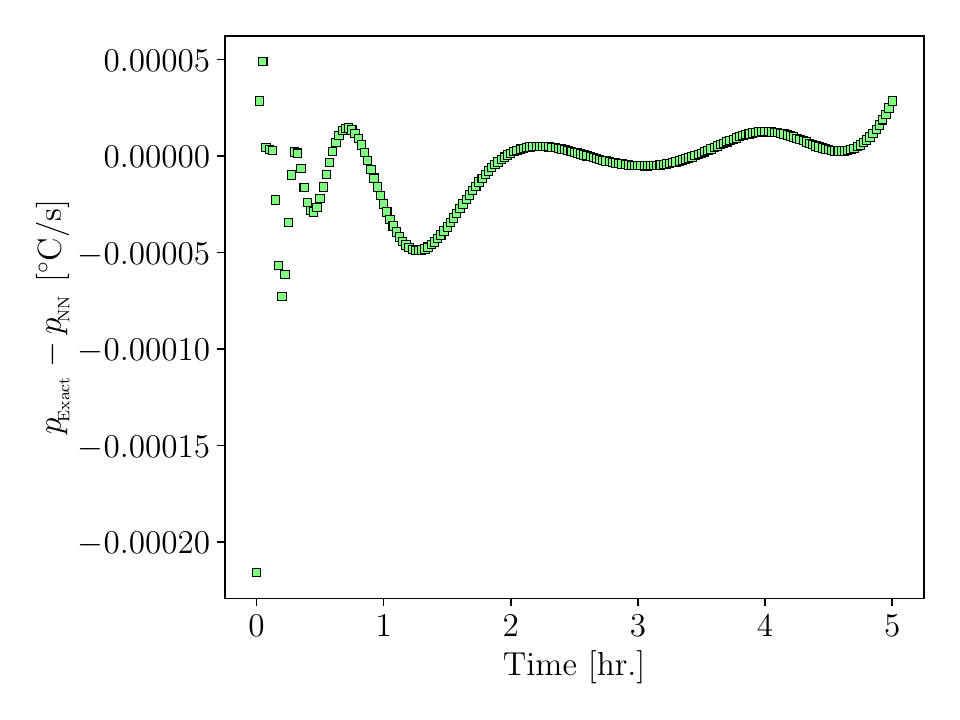}
		\caption{Difference in estimated pseudoforce}\label{fig:solar_p_error}
	\end{subfigure}\\
	\begin{subfigure}[t]{0.45\textwidth}
		\centering
		\includegraphics[scale=0.45]{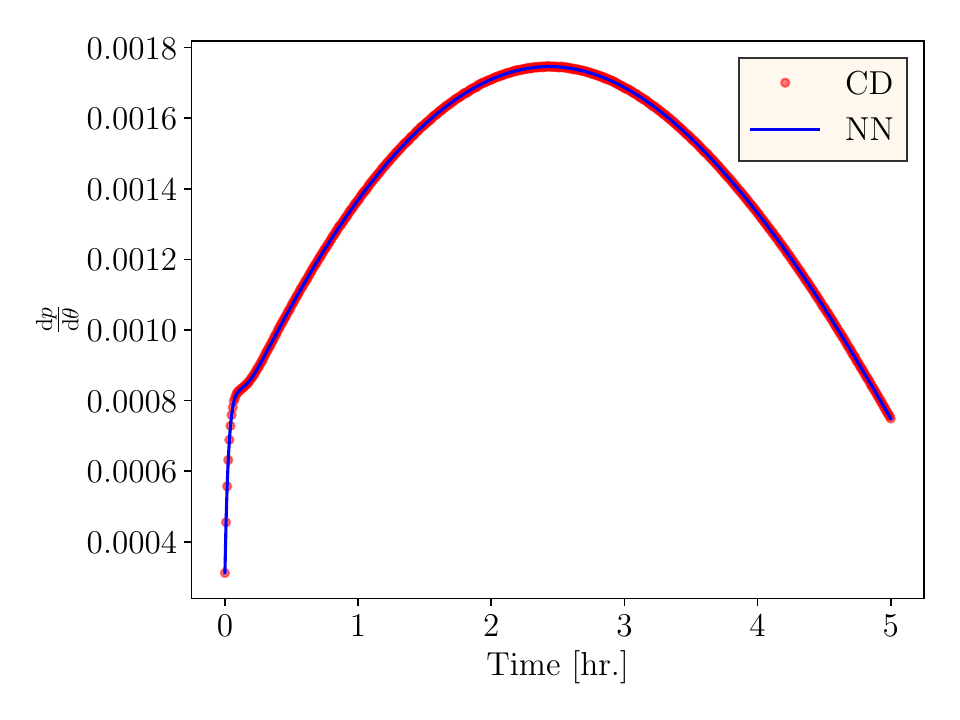}
		\caption{Comparison of the gradient of the pseudoforce}\label{fig:solar_dpdtheta_comp}
	\end{subfigure}
	\hfill
	\begin{subfigure}[t]{0.45\textwidth}
		\centering
		\includegraphics[scale=0.45]{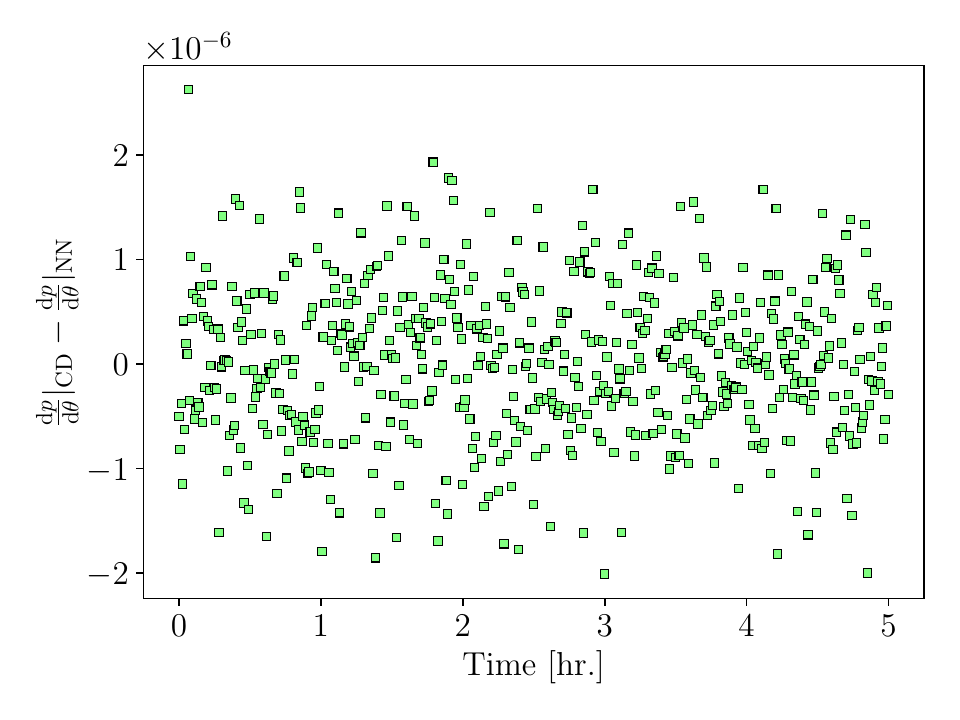}
		\caption{Difference in the gradient of the pseudoforce}\label{fig:solar_dpdtheta_error}
	\end{subfigure}
    \caption{Comparison of the neural network predictions of the pseudoforce and its gradient for the nonlinear thermal model with the solution of \eqref{eq:nvie_sens} and central difference in Example I. 
    Here ``NN'' refers to quantities estimated using neural network surrogates; ``Exact'' refers to the solution of \eqref{eq:nvie}; and ``CD'' denotes gradient estimated using central difference with a spacing of $10^{-3}$ $m^2$ for $A_c$. 
    }
	\label{fig:solar_p} 
\end{figure} 
Figure \ref{fig:solar_p} compares the neural network predictions of the nonlinear pseudoforce and its corresponding design sensitivity with the reference solutions obtained from the nonlinear governing equations. As shown in Figures \ref{fig:solar_p_comp} and \ref{fig:solar_p_error}, the neural network accurately reproduces the nonlinear pseudoforce throughout the entire five hour operating period.
Figure \ref{fig:solar_dpdtheta_comp} compares the gradient of the pseudoforce $\mathrm{d}p/\mathrm{d}\theta$ obtained through automatic differentiation with the reference gradient computed using the central difference method. The two curves exhibit excellent agreement over the entire operating time. 
The corresponding gradient error, which is on the order of $10^{-6}$, is shown in Figure \ref{fig:solar_dpdtheta_error}. 

\begin{table}[!htb]
\centering
\caption{Optimization parameters used in Example I.}
\label{tab:opt_settings}
\begin{tabular}{lc}
\hline \Tstrut
Parameter & Value \Bstrut\\
\hline \Tstrut
Design bounds & $15 \leq A_c \leq 45$ m$^2$ \\
$T_s^{\min}$ & $29.5^\circ$C \\
$T_s^{\max}$ & $90^\circ$C \\
$T_c^{\max}$& $75^\circ$C \\
$t_{\max}$ & 5 hrs.\Bstrut\\ 
\hline
\end{tabular}
\end{table}

\begin{table}[!htb]
\centering
\caption{Comparison of optimization performance for different Monte Carlo sample sizes and variance weights.
}
\label{tab:MCcomparison}
























\begin{tabular}{ccccc}
\hline
\Tstrut
Case & MC Samples, $n_s$ & Optimal $A_c$ (m$^2$) & Dynamics Eval. & Time (s) \Bstrut\\
\hline

\Tstrut
\multirow{7}{2.7 cm}{$\lambda_J=\lambda_{h_j}=0$, $j=1,2,3$}
& 10
& 36.521
& $4.2\times10^2$
& 6.26
\\

& 25
& 36.488
& $8.0\times10^2$
& 11.64
\\

& 50
& 37.195
& $2.0\times10^3$
& 28.92
\\

& 100
& 38.225
& $4.0\times10^3$
& 58.68
\\

& 500
& 40.419
& $1.9\times10^4$
& 276.96
\\

& 2,000
& 39.962
& $8.8\times10^4$
& 1215.76
\\

& 5,000
& 40.302
& $1.9\times10^5$
& 2597.21
\Bstrut\\

\hline

\Tstrut
\multirow{7}{2.7cm}{$\lambda_J=\lambda_{h_j}=5$, $j=1,2,3$}
& 10
& 32.983
& $3.8\times10^2$
& 5.91
\\

& 25
& 34.926
& $6.5\times10^2$
& 11.97
\\

& 50
& 34.742
& $1.9\times10^3$
& 28.21
\\

& 100
& 35.287
& $4.6\times10^3$
& 67.31
\\

& 500
& 35.676
& $2.3\times10^4$
& 332.23
\\

& 2,000
& 36.085
& $8.0\times10^4$
& 1130.30
\\

& 5,000
& 36.297
& $3.2\times10^5$
& 4825.38
\Bstrut\\

\hline

\Tstrut
\multirow{7}{2.7 cm}{$\lambda_J=100,\;\lambda_{h_j}=5$, $j=1,2,3$}
& 10
& 16.802
& $6.2\times10^2$
& 5.62
\\

& 25
& 17.298
& $2.2\times10^3$
& 22.88
\\

& 50
& 16.380
& $2.8\times10^3$
& 22.19
\\

& 100
& 16.851
& $6.6\times10^3$
& 57.96
\\

& 500
& 16.842
& $3.0\times10^4$
& 439.97
\\

& 2,000
& 15.602
& $8.8\times10^4$
& 1235.86
\\

& 5,000
& 15.604
& $1.8\times10^5$
& 2508.74
\Bstrut\\

\hline
\end{tabular}
\end{table}
The trained neural network surrogate for $p(t;\thetaa,\xii)$ in \eqref{eq:solar_g} is used with the gradient-based optimization framework that uses the SLSQP algorithm to solve the optimization problem \eqref{eq:op_ex1} to estimate the optimized collector area using the values reported in Table \ref{tab:opt_settings}. To investigate the influence of uncertainty sampling on the optimization process, the proposed framework is evaluated using different numbers of Monte Carlo samples per iteration, ranging from $10$ to $5,000$ samples.

\begin{figure}[!htb]
    \centering
    \includegraphics[width=0.5\linewidth]{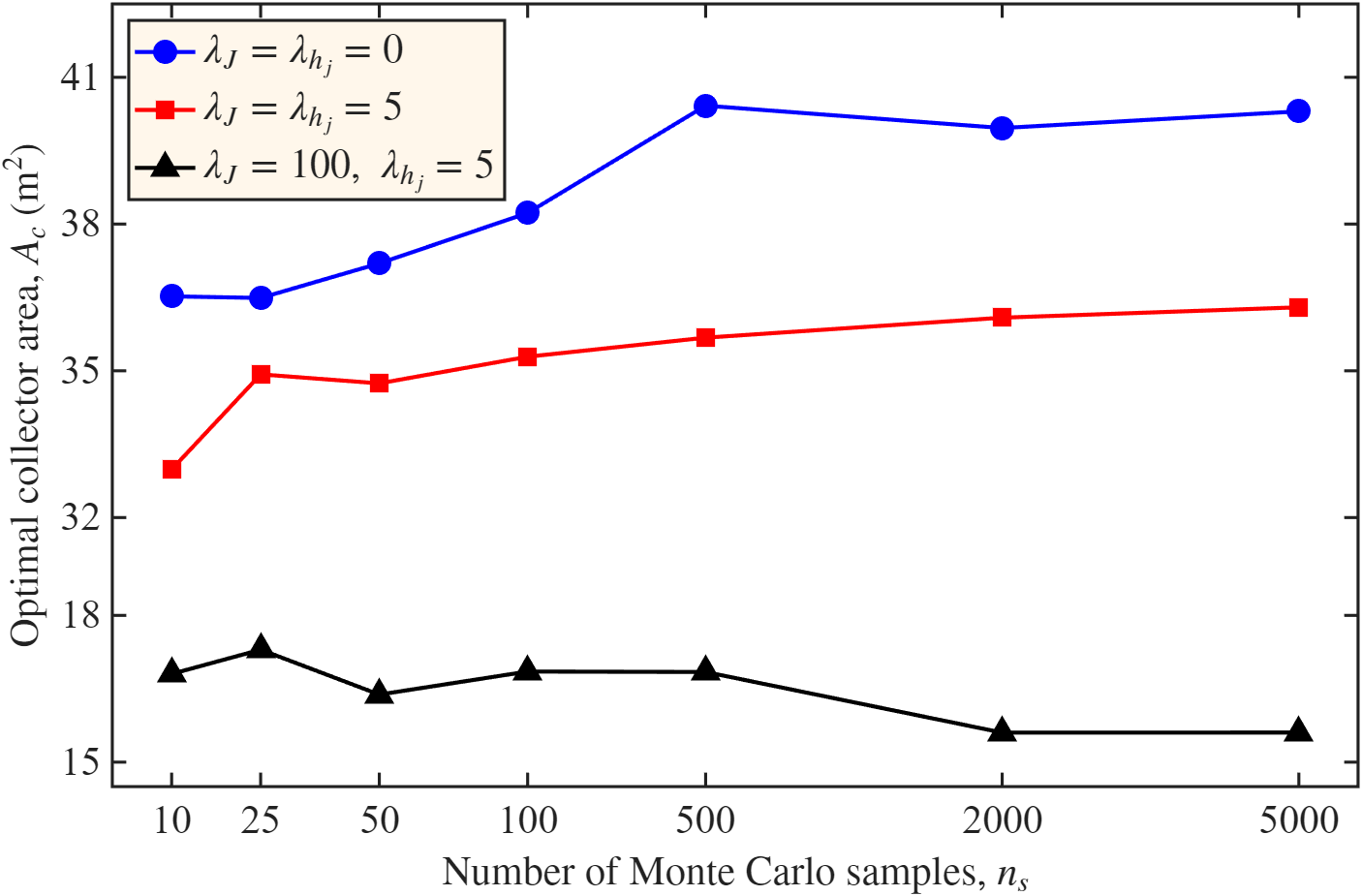}
    \caption{Effect of Monte Carlo sample size on the optimal design in Example I.}
    \label{fig:AC_MC}
\end{figure}

Table~\ref{tab:MCcomparison} and Figure~\ref{fig:AC_MC} compare the optimization results obtained using different Monte Carlo sample sizes and variance weights. The results show that the optimized collector area varies significantly when a small number of MC samples is used, whereas the solutions obtained with 2,000 and 5,000 i.i.d. samples per iteration are much closer. However, the computational cost increases considerably with the number of Monte Carlo samples because each realization requires a neural network evaluation and state space simulation. 
This suggests that 2,000 samples provide a practical compromise between computational cost and empirical solution stability. 
For the case $\lambda_J=\lambda_{h_j}=0$ (\textit{i.e.}, with no contributions from the variance), the optimal collector area varies with the number of MC samples, reaching $40.419$ at $500$ samples and decreasing to $39.962$ at 2,000 samples. Increasing the sample size further to 5,000 results in only a small change to $40.302$, indicating that the optimal collector area stabilized at approximately 2,000 MC samples. When variance is included in both the objective and constraints, with $\lambda_J=\lambda_{h_j}=5$, the corresponding optimal collector areas are $35.676$, $36.085$, and $36.297$~m$^2$, respectively. A much stronger objective variance penalty, $\lambda_J=100$ with $\lambda_{h_j}=5$, results in substantially smaller optimal designs of $16.842$, $15.602$, and $15.604$~m$^2$ for $500$, 2,000, and 5,000 samples, respectively. The very small difference between the 2,000 and 5,000 sample solutions in this case further indicates convergence with respect to the number of samples. Overall, these results demonstrate the trade-off between maximizing the expected thermal performance and reducing its sensitivity to uncertainty, with larger variance weights favoring smaller and more robust collector designs. 

Once trained, the neural network requires 20.26 and 20.59 minutes to complete the optimization for the first and third cases, respectively, using 2,000 MC samples. In comparison, using a standard nonlinear solver such as SciPy's \texttt{RK45} to evaluate the solar system is projected to require 10.85 hours for the $8.8\times10^4$ model evaluations needed during optimization. For the second case, the proposed approach requires 18.84 minutes, compared with a projected 9.86 hours for the $8.0\times10^4$ solar system evaluations using SciPy's \texttt{RK45}. These results show that, once the neural network is trained, the proposed approach substantially reduces the computational cost of robust optimization for the nonlinear solar system considered in this example. 
 The reported computational times are obtained on a laptop equipped with an AMD Ryzen 7 8845HS processor (3.80 GHz), 16 GB of RAM, running 64-bit Windows 11.

\begin{figure}[!htb]
	\centering
	
	\begin{subfigure}[t]{0.48\textwidth}
		\centering
		\includegraphics[width=\textwidth]{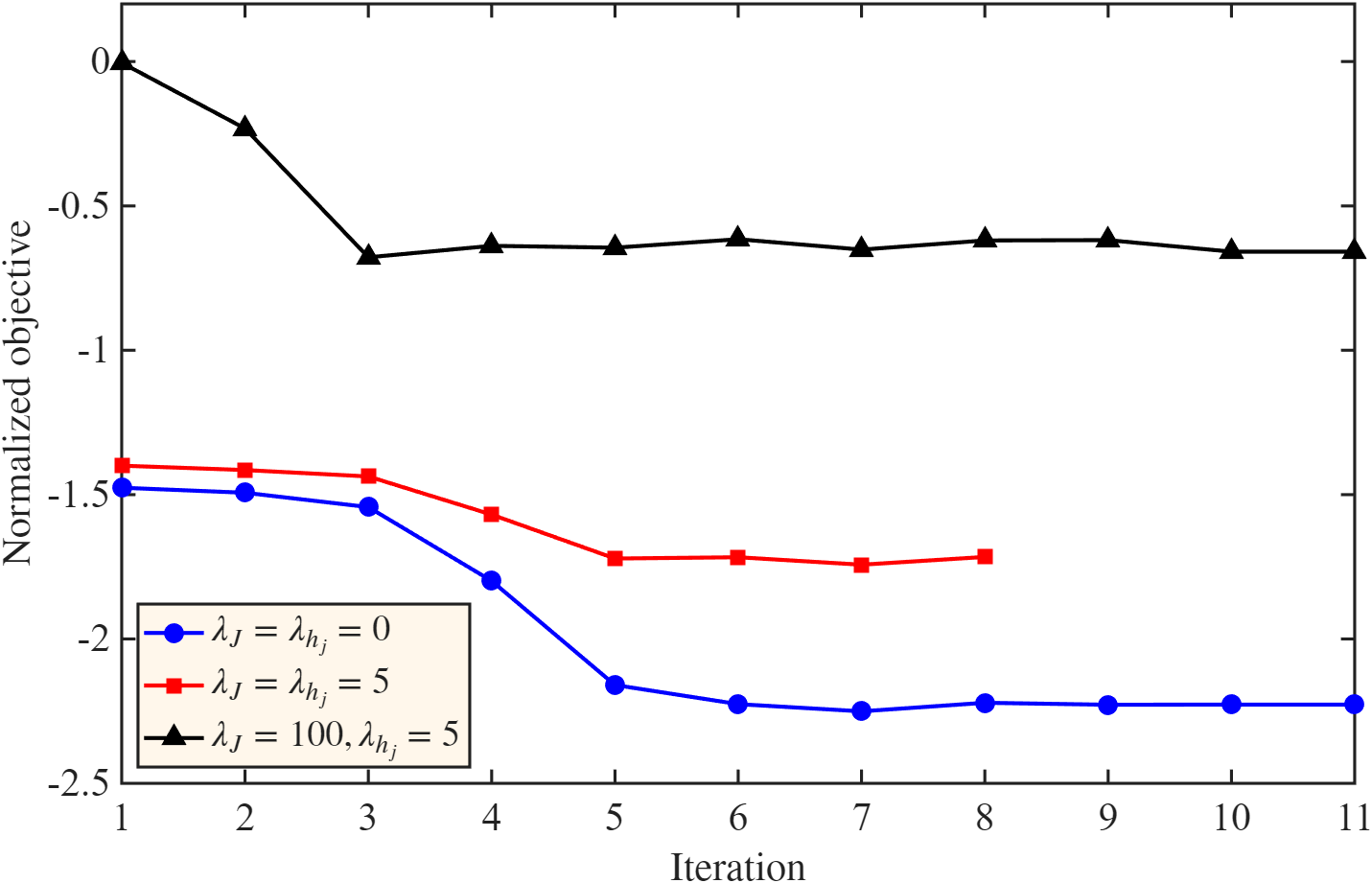}
		\caption{Normalized objective function.    
        }
		\label{fig:obj_con_a}
	\end{subfigure}
	\hfill
	\begin{subfigure}[t]{0.47\textwidth}
		\centering
		\includegraphics[width=\textwidth]{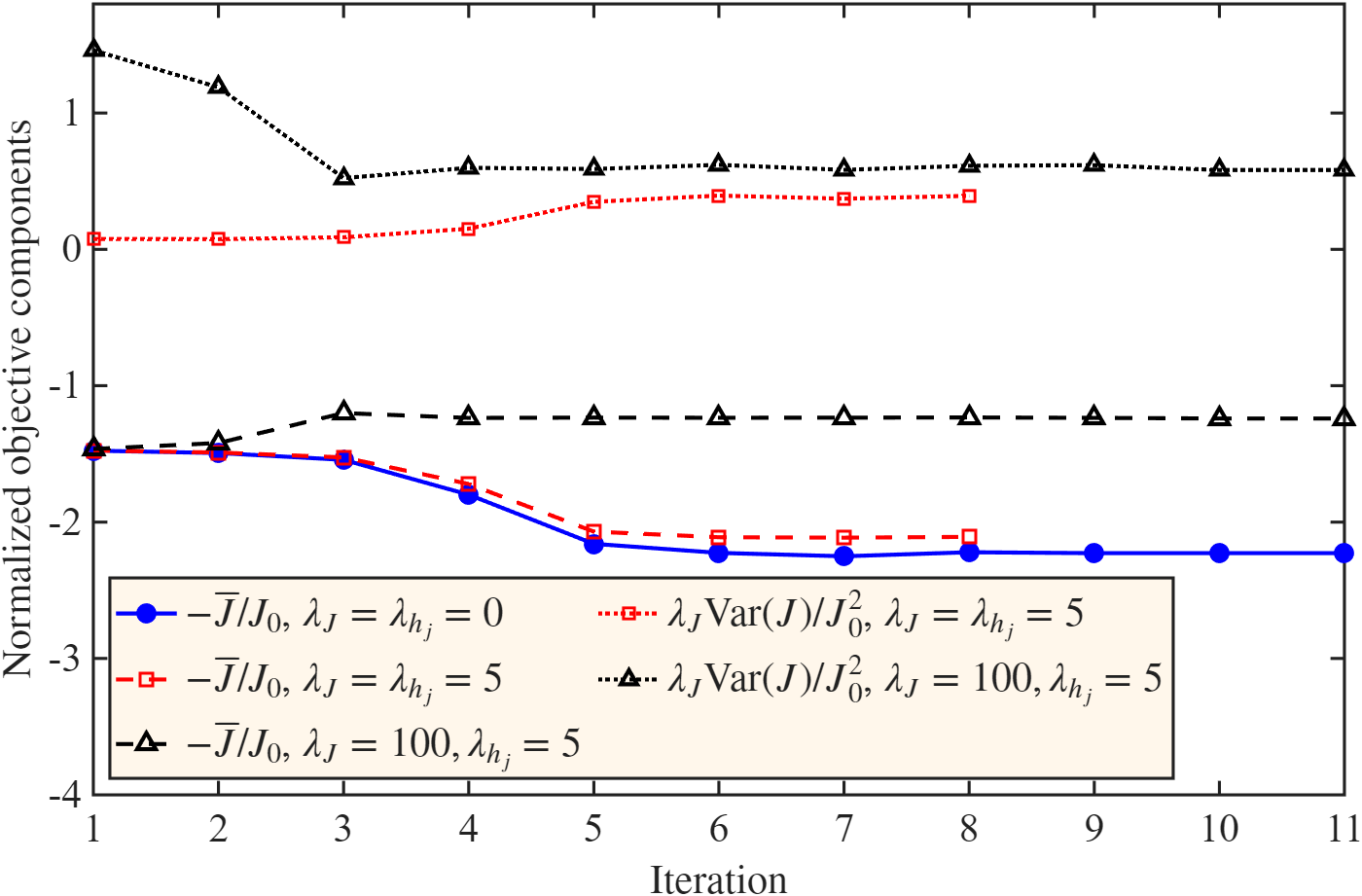}
		\caption{Normalized mean objective and weighted variance contribution. 
        }
		\label{fig:obj_con_b}
	\end{subfigure}
	
	\caption{Convergence history of the thermal system optimization 
	using 2,000 Monte Carlo samples per iteration.}
	\label{fig:obj_con}
	
\end{figure}

The convergence history of the objective for $n_s=$~2,000 uncertainty samples is presented in Figure \ref{fig:obj_con}. In particular, Figure \ref{fig:obj_con}(a) shows the convergence of the total normalized objective function for the three weighting cases, while Figure \ref{fig:obj_con}(b) presents the corresponding mean and weighted variance contributions. The total objective decreases and converges in all three cases. However, the relative contributions of the mean and variance depend on the selected variance weight, $\lambda_J$. When $\lambda_J$ is small, the optimizer places greater emphasis on reducing the mean objective, while allowing a larger variance. In contrast, increasing $\lambda_J$ places greater emphasis on reducing the variability of the objective. Consequently, for the case with the largest $\lambda_J$, the optimizer accepts a higher mean value in exchange for a smaller variance contribution. These results demonstrate the trade-off between minimizing the expected objective and reducing its sensitivity to uncertainty.

\begin{figure}[!htb]
    \centering

    \begin{subfigure}[t]{0.48\textwidth}
        \centering
        \includegraphics[width=\textwidth]{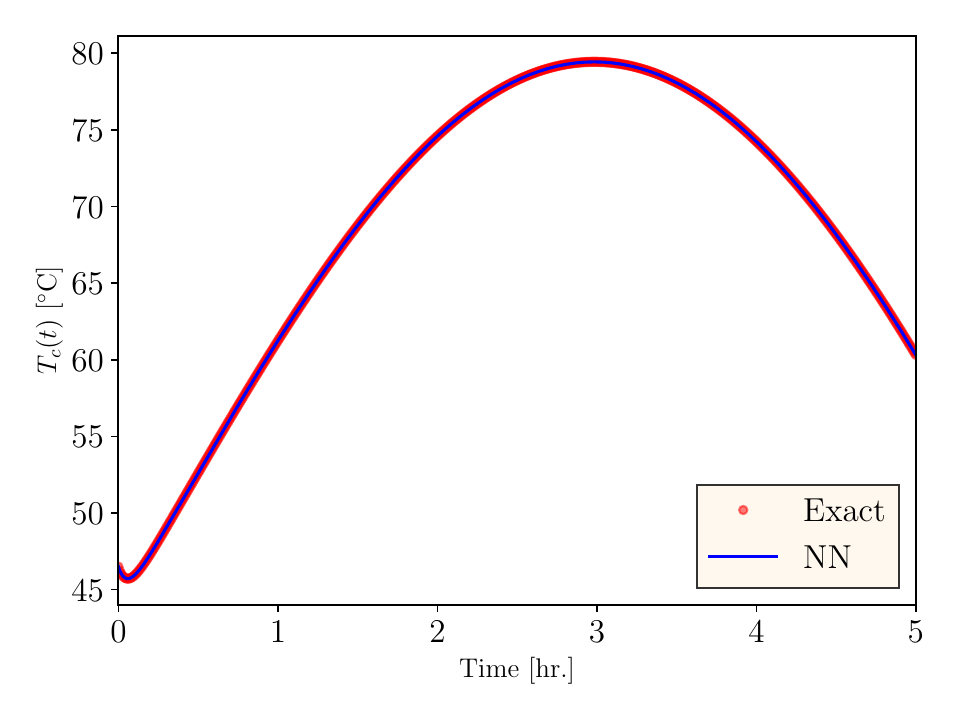}
        \caption{Collector temperature, $T_c(t)$.}
        \label{fig:solar_Tc}
    \end{subfigure}
    \hfill
    \begin{subfigure}[t]{0.48\textwidth}
        \centering
        \includegraphics[width=\textwidth]{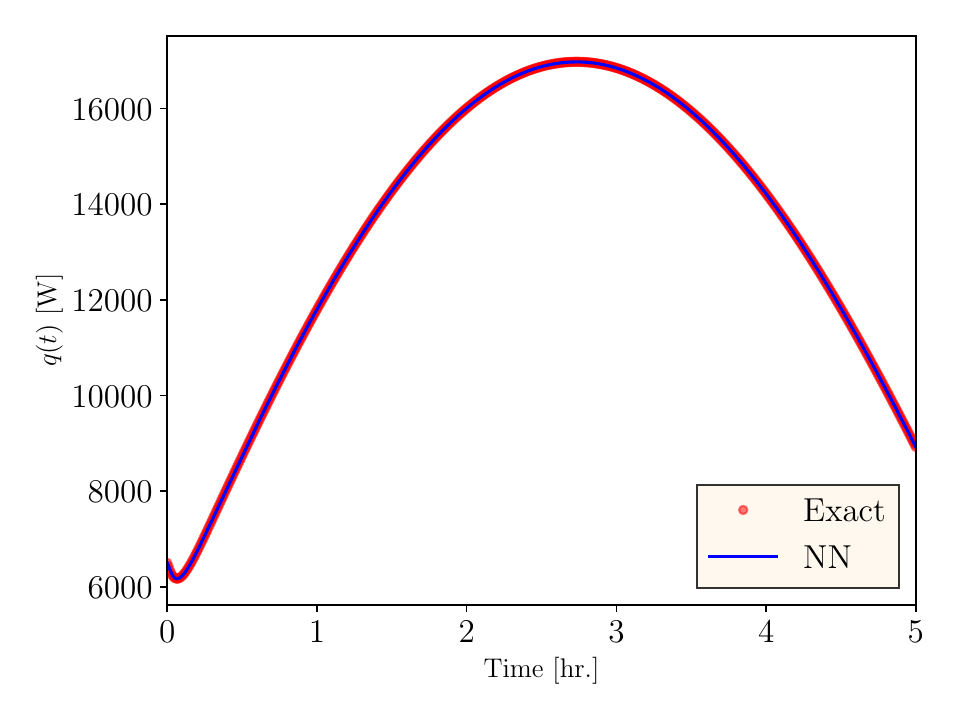}
        \caption{Heat transfer rate, $q(t)$.}
        \label{fig:solar_q}
    \end{subfigure}

    \vspace{0.3cm}

    \begin{subfigure}[t]{0.48\textwidth}
        \centering
        \includegraphics[width=\textwidth]{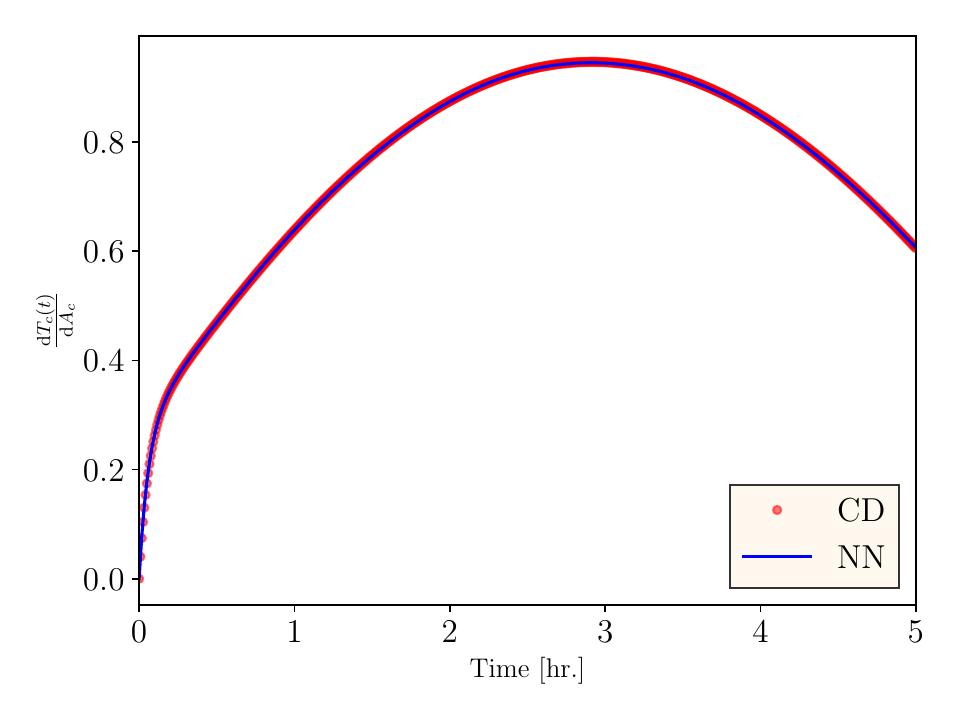}
        \caption{Temperature sensitivity.}
        \label{fig:solar_dTc}
    \end{subfigure}
    \hfill
    \begin{subfigure}[t]{0.48\textwidth}
        \centering
        \includegraphics[width=\textwidth]{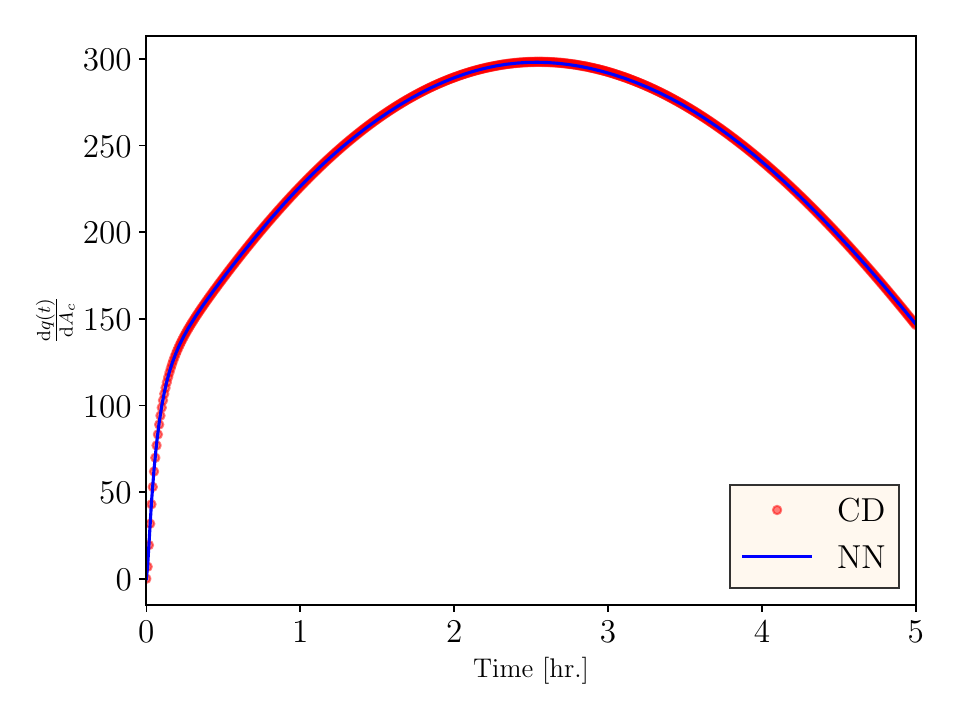}
        \caption{Heat transfer rate sensitivity.}
        \label{fig:solar_dq}
    \end{subfigure}

    \caption{Comparison of the proposed surrogate model and the
    high-fidelity nonlinear solution for Example I at the optimized
    collector area. The time-dependent responses and their
    sensitivities with respect to the collector area are evaluated
    for one realization of the uncertain parameters.
    }

    \label{fig:solar_validation}
\end{figure}

To evaluate the accuracy of the proposed surrogate at the optimized design, Table \ref{tab:opt_verification} compares the mean and variance components of the objective and mean active constraint values obtained using the proposed bi-fidelity method with the high-fidelity reference computed by solving the nonlinear high-fidelity model. 
At the optimized collector area, the proposed approach and high-fidelity models yield objective mean components of \(-2.1063\) and \(-2.1062\), objective variance components of \(0.0711\) and \(0.0710\), and active constraint values of \(-0.0038\) and \(-0.0039\), respectively, as reported in Table~\ref{tab:opt_verification}. 
Thus, the surrogate preserves both the objective statistics and the active constraint at the optimized design. 
To further examine the accuracy of the surrogate model, the
time dependent collector temperature, $T_c(t)$, and heat transfer
rate, $q(t) = \rho_s c_s \upsilon_s\left(T_{hc}(t;\theta,\xii)-T_s(t;\theta,\boldsymbol{\xi})\right)$, together with their sensitivities with respect to the collector area, $A_c$, are evaluated for one realization at the optimized design. As shown in
Figure~\ref{fig:solar_validation}, the responses and gradients predicted by the proposed method closely follow those obtained from the high-fidelity nonlinear model. 

\begin{table}[!htb]
\centering
\caption{
Comparison of the proposed bi-fidelity method with the high-fidelity (HF) reference at the optimized designs using $n_s=2,000$ Monte Carlo samples for Example~I and Example~II. 
} 
\label{tab:opt_verification}
\begin{tabular}{llccc}
\hline \Tstrut
Example & Quantity & Proposed & HF reference 
\Bstrut\\
\hline \Tstrut
 & Mean component of objective & $-2.1063$ & $-2.1062$ 
\\
I ($\lambda_J=\lambda_{hj}=5$) &  Variance component of objective & ~~~0.0711 & ~~0.0710 
\\
 & Mean active constraint value ($h_3$) & $-0.0038$ & $-0.0039$ 
 \Bstrut\\ \hline 
\Tstrut & Mean component of objective & ~~~3.0649 & ~~~\!3.0863 
\\
II ($\lambda_J=0.1$, $\lambda_{h}=20$)& Variance component of objective & ~~~0.2831 & ~\!~~0.2944 
\\
 & Mean base displacement constraint & $-0.0036$ & $-0.0037$ 
 \Bstrut\\ 
\hline
\end{tabular}
\end{table}

\subsection{Example II: 11-story 2-bay Base-isolated Structure} 
\subsubsection{Description}
An 11-story 2-bay structure with a hysteretic base isolation layer adapted from Kamalzare et al. \cite{kamalzare2015efficient} and De et al. \cite{de2018computationally} is used in the second numerical example in this paper. 
Consistent mass matrices are used for the beams while ignoring the weights of the columns. The base layer moves horizontally with rigid in-plane behavior. With these assumptions, the governing differential equations of the base-isolated structure are given by
\begin{equation}\label{eq:11story}
\begin{split}
&\Mm \ddot{\um}_\mathrm{s} + \Cms \dot{\um}_\mathrm{s} + \Km \um_\mathrm{s} = - \Mm \rr \ddot{u}_\mathrm{g} + \Cms \rr \dot{u}_\mathrm{b} + \Km \rr u_\mathrm{b};\\
&m_\mathrm{b}\ddot{u}_\mathrm{b}+\left(\rr^T\Cms\rr\right)\dot{u}_\mathrm{b}+\left(\rr^T\Km\rr\right)u_\mathrm{b}+f_\mathrm{b}=-m_\mathrm{b}\ddot{u}_\mathrm{g}+\rr^T\Cms\dot{\um}_\mathrm{s}+\rr^T\Km\um_\mathrm{s}, 
\end{split}
\end{equation}
where $\Mm$ is the mass matrix, $\Cms$ is the damping matrix, and $\Km$ is stiffness matrix of the superstructure; $\um_\mathrm{s}$ is displacement of the superstructure relative to the ground; $m_\mathrm{b}$ is the mass of the base layer; 
$f_\mathrm{b}$ is the restoring force in the hysteretic base layer; $\ddot{u}_\mathrm{g}$ is the ground acceleration; and the influence vector for the ground acceleration is $\rr = [1,0,0,\dots,1,0,0]^T$, where the ones correspond to the horizontal displacement DOF. 
For the superstructure Rayleigh damping with 3\% damping ratios for the 1$^\mathrm{st}$ and 10$^\mathrm{th}$ is assumed. 
The structure has in total 100 DOF with three DOF per node for 33 nodes of the superstructure and one horizontal DOF for the base layer. 
The ground acceleration used here is a record of the $1940$ El Centro earthquake measured at the N-S Imperial Valley Irrigation District substation with peak ground acceleration $0.348g$ is used as the ground acceleration $\ddot{u}_\mathrm{g}$. 

\begin{figure}[!htb]
	\centering
	\begin{subfigure}[t]{0.45\textwidth}
		\centering
		\includegraphics[scale=0.5]{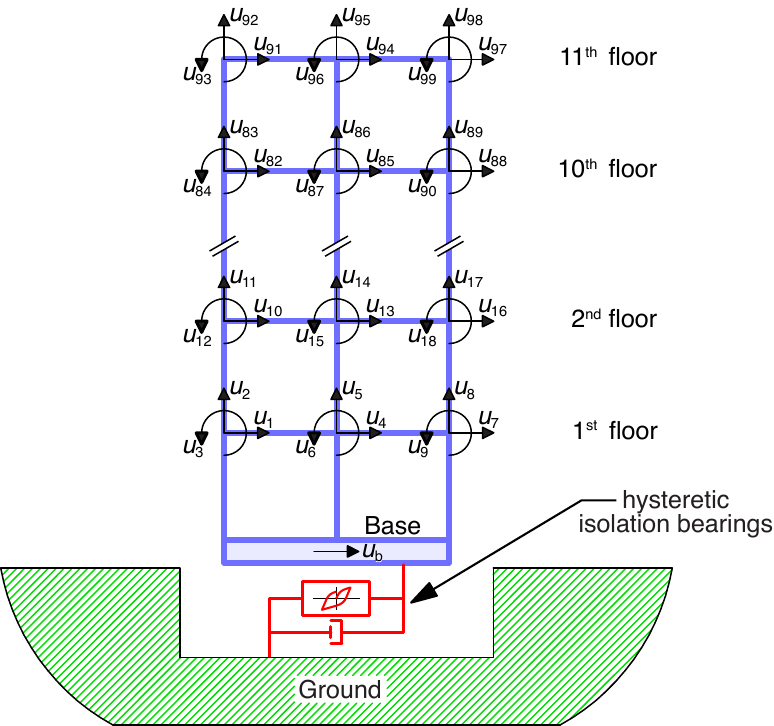}
		\caption{A 100 degree of freedom structure with a hysteretic base isolation layer.}\label{fig:100dofa}
	\end{subfigure}
	\hfill
	\begin{subfigure}[t]{0.5\textwidth}
		\centering
		\includegraphics[scale=1]{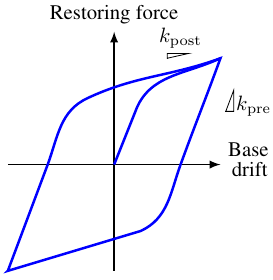}
		\caption{Bouc-Wen model for hysteresis.}\label{fig:100dofb}
	\end{subfigure}
	\caption{ The hysteretic base isolation layer described using the Bouc-Wen hysteresis model of an 11-story 2-bay 100 degree of freedom structure (adapted from Kamalzare et al. \cite{kamalzare2015efficient} and De et al. \cite{de2018computationally}) is designed in Example II.}
	\label{fig:100dof}
\end{figure}
The restoring force in the hysteretic base layer is described using the Bouc-Wen model (see Figure \ref{fig:100dofb}) \cite{wen1976method} given by
\begin{equation} 
f_\mathrm{b} = c_\mathrm{b}\dot{u}_\mathrm{b}+\kpost{u}_\mathrm{b}+\alpha z_\mathrm{hyst},
\end{equation}
where the base layer damping coefficient is $c_\mathrm{b}$; $\kpre$ and $\kpost$ are pre- and post-yield stiffness, respectively; $\alpha = \Qy(1-r_k)$ is the peak of the nonelastic force; $\Qy$ is the yield force of the hysteretic layer; $r_k=\kpost/\kpre$ is the hardening ratio; and $z_\mathrm{hyst}$ is an auxiliary variable that describes the evolution of the hysteretic loop. The evolution of $z_\mathrm{hyst}$ is given by
\begin{equation} \label{eq:zhyst}
\dot{z}_\mathrm{hyst} = A \dot{u}_\mathrm{b} - \beta \dot{u}_\mathrm{b} |z_\mathrm{hyst}| - \gamma z_\mathrm{hyst} |\dot{u}_\mathrm{b}|. 
\end{equation}
In \eqref{eq:zhyst}, the parameters are chosen as $A=2\beta=2\gamma=\kpre/\Qy$, which produces identical loading and unloading curves \cite{ramallo2002smart} with $z_\mathrm{hyst}$ in $[-1,1]$. 
To construct the nominal linear system, the restoring force is decomposed as
\begin{equation}
f_\mathrm{b}=\bar{c}\dot{u}_\mathrm{b}+\bar{k}u_\mathrm{b}+p(t;\thetaa,\xii),
\end{equation}
where $\bar{c}=35~\mathrm{kN\,s/m}$ and $\bar{k}=750~\mathrm{kN/m}$ are the nominal damping and stiffness values, respectively. Consequently, the nonlinear and uncertain pseudoforce is given by
\begin{equation}
p(t;\thetaa,\xii)=(c_\mathrm{b}-\bar{c})\dot{u}_\mathrm{b}(t)+(\kpost-\bar{k})u_\mathrm{b}(t)+\alpha z_\mathrm{hyst}(t). 
\label{eq:building_pseudoforce}
\end{equation} 
The governing equations of this structure are then converted to the state-space form in \eqref{eq:state-space} using $\Xm=[\um_\mathrm{s}^T\quad u_\mathrm{b}\quad\dot{\um}_\mathrm{s}^T\quad\dot{u}_\mathrm{b}]^T$, with $\bar{c}$ and $\bar{k}$ incorporated into the nominal linear system and $p(t;\thetaa,\xii)$ treated as the nonlinear and uncertain pseudoforce. 
Three parameters, namely, $c_\mathrm{b}$, $\kpost$, and $r_k$ of the restoring isolation force are assumed uncertain with their probability distributions given in Table \ref{tab:ExI_param} following De et al. \cite{de2018computationally}. These distributions are chosen such that significant variations are observed in the response. Further, lognormal distribution is assumed for $\kpost$ and Gaussian distribution truncated below at zero for $c_\mathrm{b}$ to ensure they are positive. 

\subsubsection{Optimization Problem} 

The yield force of the isolation layer $\Qy$ in \% of the total weight of the structure is used as the design parameter in this example for optimization (\textit{i.e.}, $\theta=\Qy$). The lower and upper bounds of this design parameter are assumed as 2.0\% and 8.0\% of the total weight of the structure (see Table \ref{tab:ExI_param}). 
The aim of this optimization is to reduce the horizontal roof acceleration of the structure under ground excitation. Hence,
for each realization of the uncertain parameters, the cost function is defined as
\begin{equation}
    {J}(\theta,\boldsymbol{\xi})
    =
    \frac{1}{2}
    \left\|
    Y(t;\theta,\boldsymbol{\xi})
    \right\|_2^2
    =
    \frac{1}{2}
    \sum_{i=1}^{n_t}
    \left[
    \ddot{u}_{97}(t_i;\theta,\boldsymbol{\xi})
    \right]^2,
\end{equation}
where $n_t$ is the number of time steps used in discretizing the time domain. 
The maximum absolute horizontal displacement of the base layer is used to define the feasibility condition as
\begin{equation}
    {h}(\theta,\boldsymbol{\xi})
    =
    \max_{t}
    \left|
    u_b(t;\theta,\boldsymbol{\xi})
    \right|
    -
    u_{b,\max}
    \leq 0,
\end{equation} 
where $u_{b,\max}=0.1$~m is the allowable maximum horizontal base displacement. 
To improve the numerical scaling of the optimization problem, the objective function is normalized by a reference value $J_0=100$ $\mathrm{m}^2/\mathrm{s}^4$ and the constraint is scaled using $h_0=u_{b,\max}$. 
The normalized robust optimization problem is then formulated as
\begin{equation}
    \begin{aligned}
    \min_{\theta}\quad
    & 
    \mathbb{E}_{\boldsymbol{\xi}}
    \left[\frac{1}{J_0}
    J(\theta,\boldsymbol{\xi})
    \right]
    +
    \lambda_J
    \operatorname{Var}
    \left[\frac{1}{J_0}
    J(\theta,\boldsymbol{\xi})
    \right]
    \\
    \text{subject to}\quad
    &
    \mathbb{E}_{\boldsymbol{\xi}}
    \left[\frac{1}{h_0}
    h(\theta,\boldsymbol{\xi})
    \right]
    +
    \lambda_{h}
    \operatorname{Var}
    \left[\frac{1}{h_0}
    h(\theta,\boldsymbol{\xi})
    \right]
    \leq
    0,
    \\
    &
    2\% \leq \theta \leq 8\%.
    \end{aligned}
    \label{eq:op_ex2}
\end{equation} 
Note that this cost function encourages the design parameter $\Qy$ to be small. However, this increases the base displacement, which is used as a constraint to be less than $0.1$~m. We use $n_s=2,000$ similar to the previous example to provide a  practical compromise between computational cost and empirical solution stability. 
The initial value of $\Qy$ used in the optimization is set to 6.0\%.

\begin{table}
\begin{center}
	\begin{threeparttable}[!htp]
		\caption{Probability distributions of the uncertain parameters and the lower and upper bounds of the design parameter in Example II. } \label{tab:ExI_param} 
		\begin{tabular}{c c c c c c c} 
			\hline 
			\Tstrut
			 \multirow{2}{*}{Type} & Parameter & \multirow{2}{*}{Distribution} & \multirow{2}{*}{Mean} & \multirow{2}{*}{Std. Dev.} & Lower  & Upper \\ 
			  & [Unit] & & & & bound & bound \\[0.5ex] 
			\hline 
			\Tstrut
			\multirow{3}{*}{Uncertain} & $\kpost$ [kN/m] & Lognormal & 750 & 20 & 0 & $\infty$\\ 
			& $c_\mathrm{b}$ [kN$\cdot$s/m] & Truncated Gaussian$^*$ & 35 & 2.5 & 0 & $\infty$\\
			& $r_k$ & Uniform & 0.1875 & 0.0361 & 0.125 & 0.250 \Bstrut\\ \hline \Tstrut
			Design & $\Qy$ (\%)$^\dagger$ & --- & --- & --- & 2.0 & 8.0 \\ [1ex] 
			\hline 
		\end{tabular}
		\begin{tablenotes}
			\item[$^*$] Truncated below at zero; \item[$^\dagger$] in \% of the total weight of the structure.
		\end{tablenotes}
	\end{threeparttable}
\end{center}
\end{table}

\subsubsection{Results}

\begin{figure}[!htb]
	\centering
	\begin{subfigure}[t]{0.45\textwidth}
		\centering
		\includegraphics[width=\linewidth]{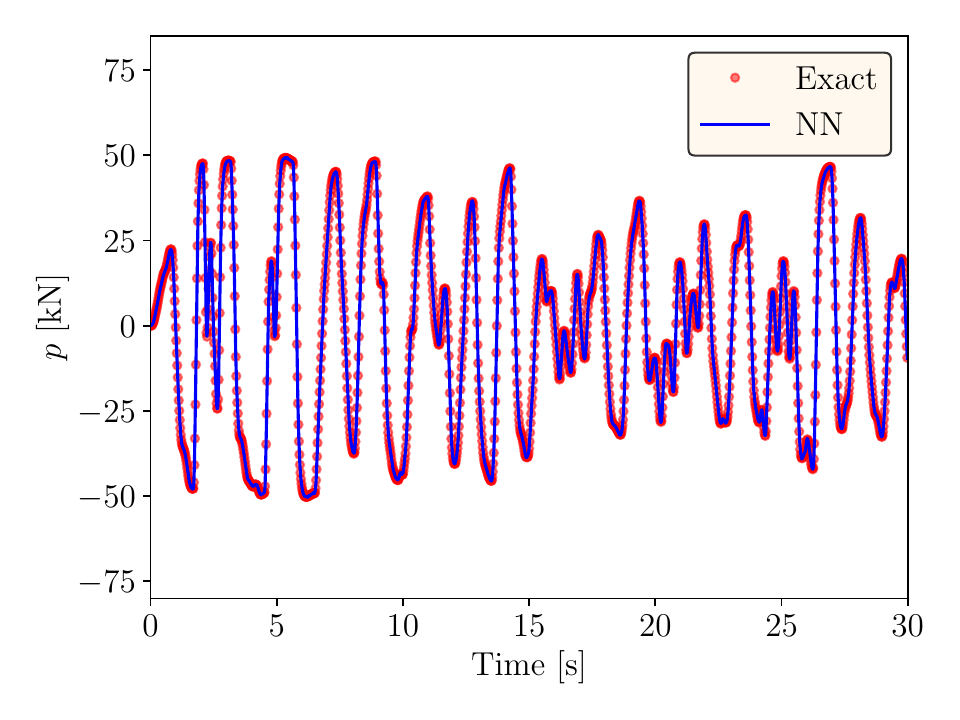}
		\caption{Comparison is estimated pseudoforce}\label{fig:100dof_p_comp}
	\end{subfigure}
	\hfill
	\begin{subfigure}[t]{0.45\textwidth}
		\centering
		\includegraphics[width=\linewidth]{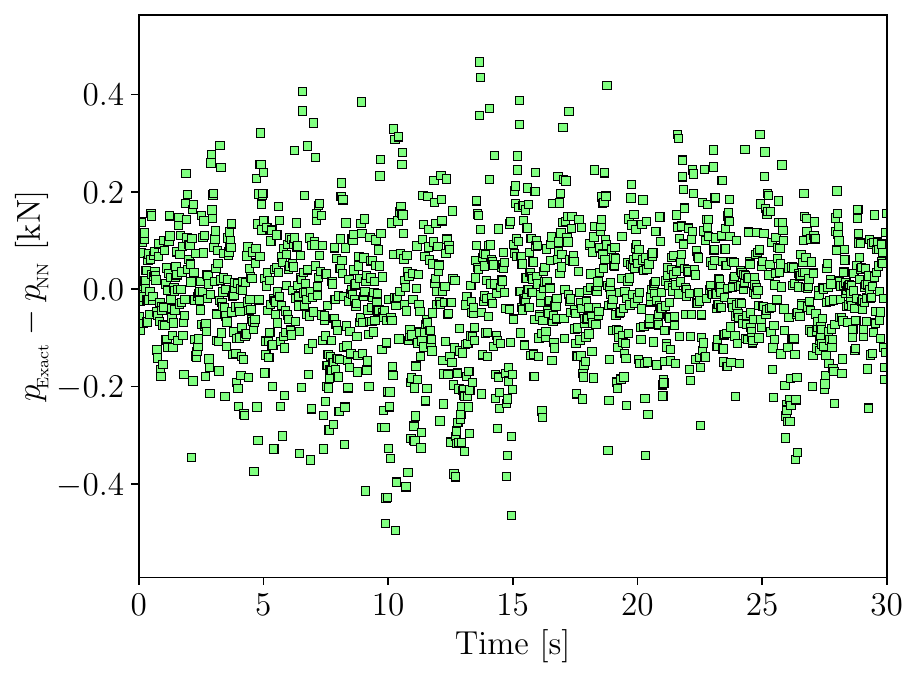}
		\caption{Difference in estimated pseudoforce}\label{fig:100dof_p_error}
	\end{subfigure}\\
	\begin{subfigure}[t]{0.45\textwidth}
		\centering
		\includegraphics[width=\linewidth]{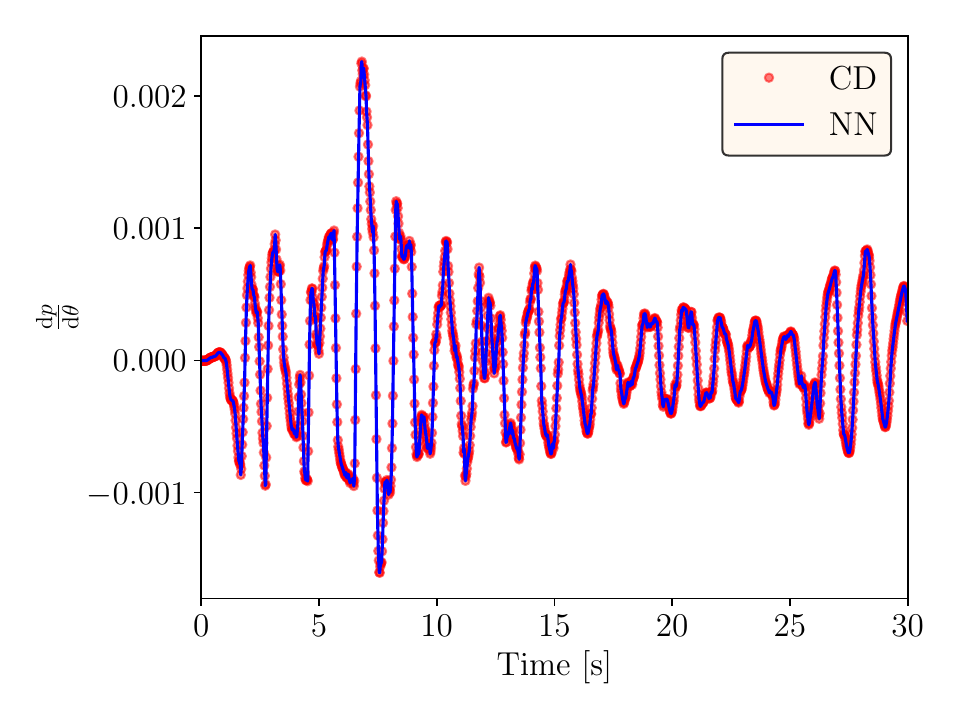}
		\caption{Comparison of the gradient of the pseudoforce}\label{fig:100dof_dpdtheta_comp}
	\end{subfigure}
	\hfill
	\begin{subfigure}[t]{0.45\textwidth}
		\centering
		\includegraphics[width=\linewidth]{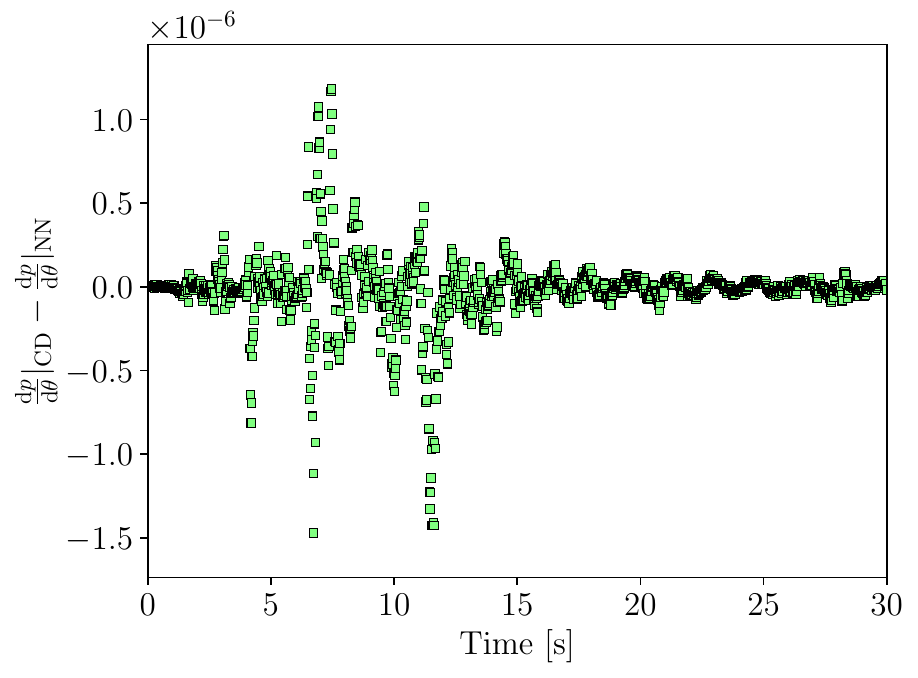}
		\caption{Difference in the gradient of the pseudoforce}\label{fig:100dof_dpdtheta_error}
	\end{subfigure}
	\caption{ Comparison of the pseudoforce and its gradients estimated using a neural network surrogate with the solution of \eqref{eq:nvie_sens} and central difference in Example II. Here ``NN" refers to quantities estimated using neural network surrogates; ``Exact" refers to the solution of \eqref{eq:nvie}; and ``CD" denotes gradient estimated using central difference with a spacing of $\Qy=0.1$\%. 
    }
	\label{fig:100dof_p}
\end{figure}

The nonlinear pseudoforce in the base layer is modeled as Equation \eqref{eq:building_pseudoforce} using a neural network.
Hence, the Bouc-Wen auxiliary variable $z_\mathrm{hyst}$ is integrated when generating the reference pseudoforce histories \cite{wojtkiewicz2011efficient,de2018computationally}. As a result, the neural network does not replace the Bouc-Wen constitutive equation directly; it approximates the resulting parameterized pseudoforce trajectory for the prescribed ground motion input.
In particular, a neural network with $N_H=4$ hidden layers with $m=100$ neurons per layer and hyperbolic tangent function 
as the activation function is used following De et al.~\cite{de2020uncertainty}. 
The network is trained using the Adam algorithm with a learning rate of $1 \times 10^{-3}$ for a maximum iteration of 10,000. This learning rate is gradually halved every 1,000 iterations. A total of 500 Monte Carlo realizations are generated to construct the dataset, of which 450 trajectories are used for training and the remaining 50 trajectories are reserved for validation and performance evaluation, with a sampling interval of $\Delta t=0.02~\mathrm{s}$. 
The training of the network takes approximately 5.3 hours. 
The relative RMSE for the validation dataset $\Dval$ using the trained network is 
$5.952\times10^{-3}$. 
Figure \ref{fig:100dof_p} compares the prediction of $p(t;\thetaa,\xii)$ and its derivative with respect to the design parameter ${\mathrm{d}p}/{\mathrm{d}\theta}$ for one realization of the uncertain parameters, where the isolator force is estimated by solving \eqref{eq:nvie} and a central difference with a spacing of 0.1\% for $\Qy$ is used for the sensitivity. This figure shows that the trained neural network accurately estimates these quantities. 

\begin{figure}[!htb]
	\centering
	
	\begin{subfigure}[t]{0.45\textwidth}
		\centering
		\includegraphics[width=\linewidth]{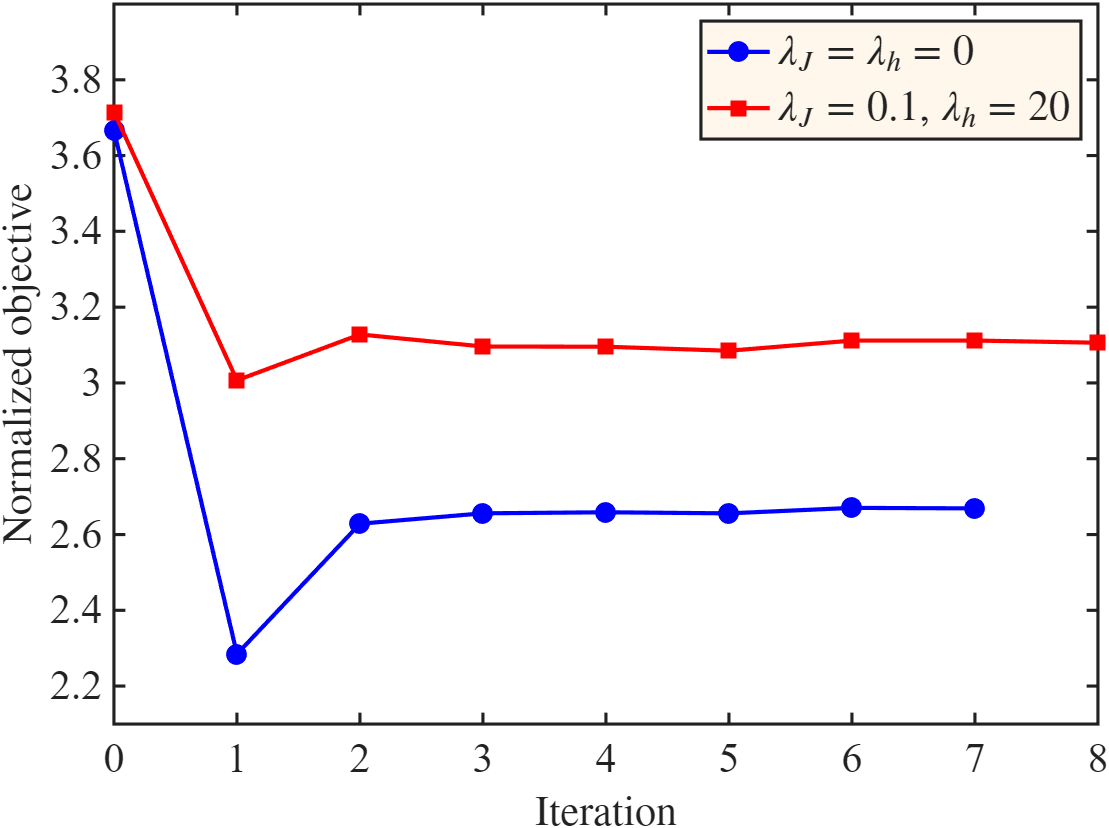}
		\caption{Normalized objective function.}
		\label{fig:structural_obj_a}
	\end{subfigure}
	\hfill
	\begin{subfigure}[t]{0.45\textwidth}
		\centering
        \includegraphics[width=\linewidth]{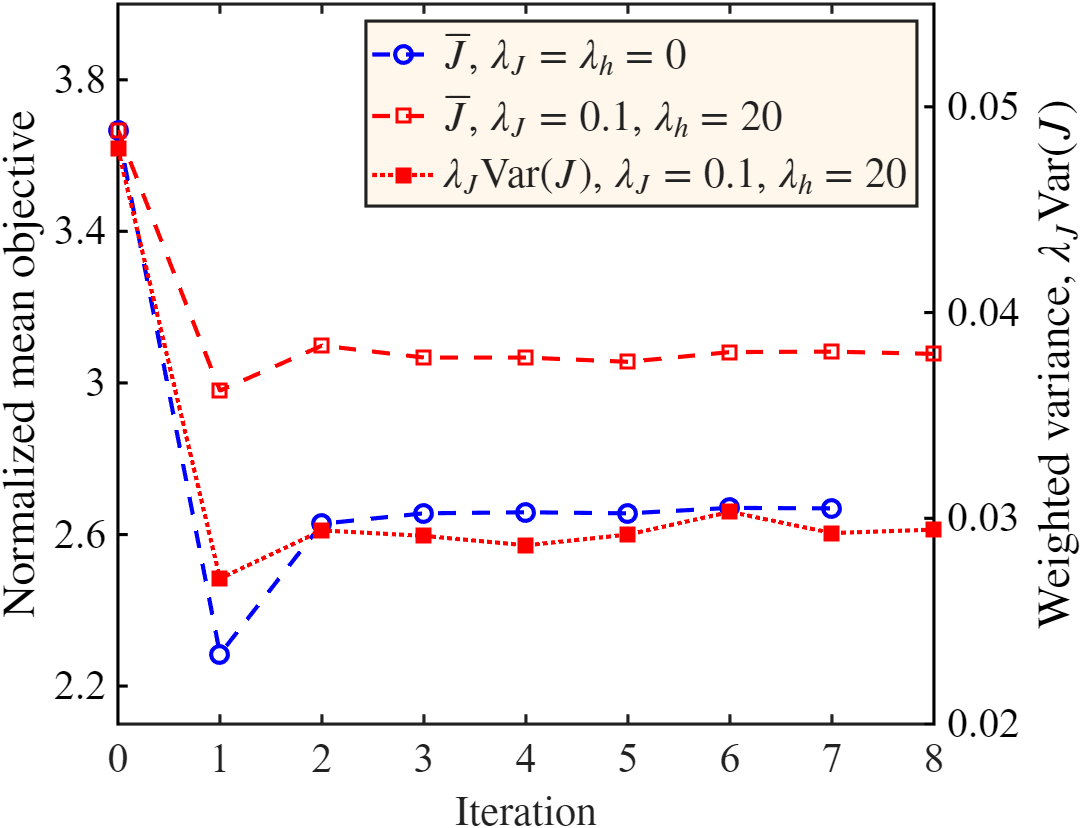}
		\caption{Normalized mean objective and weighted variance contribution.}
		\label{fig:structural_obj_b}
	\end{subfigure}
    \\
    \begin{subfigure}[t]{0.45\textwidth}
		\centering
        \includegraphics[width=\linewidth]{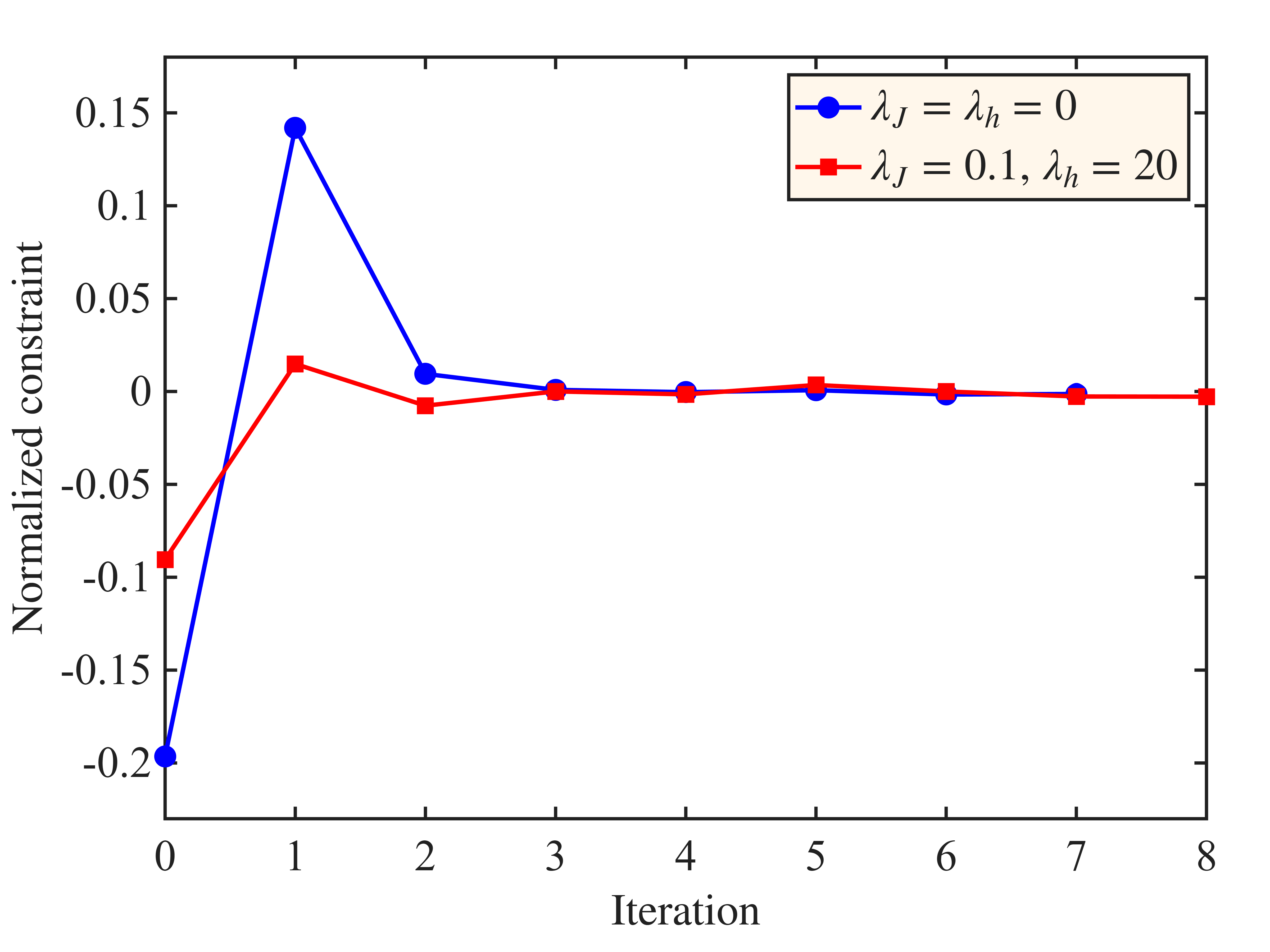}
		\caption{Normalized base displacement constraint.}
		\label{fig:structural_obj_c}
	\end{subfigure}
	
	\caption{Convergence history of the structural system optimization using 2,000 Monte Carlo samples per iteration.
    }
	\label{fig:structural_objective_convergence}
	
\end{figure}

The trained neural network is used with the gradient-based optimization framework that uses the SLSQP algorithm to solve the optimization problem \eqref{eq:op_ex2} to estimate the optimized yield force. 
Table \ref{tab:ExII_results} presents the results obtained using 2,000 Monte Carlo samples. In the first case, variance contributions are excluded from both the objective function and the constraint, \textit{i.e.}, $\lambda_J=\lambda_h=0$. The corresponding optimal yield force, $\Qy$, of the hysteretic base layer is 4.10\%. 
In the second case, $\lambda_J=0.1$ and $\lambda_{h}=20$ are used in \eqref{eq:op_ex2}. The optimal $\Qy$ in this case is 4.88 \%. 
The convergence history of the objective function for 2,000 uncertainty samples is presented in Figure~\ref{fig:structural_objective_convergence}. Figure~\ref{fig:structural_objective_convergence}(a) shows the convergence of the total normalized objective function for the two weighting cases, while Figure~\ref{fig:structural_objective_convergence}(b) presents the corresponding normalized mean objective and weighted variance contribution. Further, Figure~\ref{fig:structural_objective_convergence}(c) shows the constraint is violated initially but is satisfied at the end of the optimization. 



\begin{table}[!htb]
\begin{center}
		\caption{Optimal design parameter and optimization performance for two cases considered in Example II. 
        } \label{tab:ExII_results} 
		\begin{tabular}{c c c c c c c} 
			\hline 
			\Tstrut
			{$\lambda_J$} & {$\lambda_{h}$} & Time (hr.)  &  Dynamics Evaluation & Optimal $\Qy$ (\%)  \\ [0.5ex] 
			\hline 
			\Tstrut 
			0 & 0 & 0.97 
            & $4.4\times10^4$ & 4.10 \\
			0.1 & 20 & 1.09 
            & $4.8\times10^4$ & 4.88 \\[1ex] 
			\hline 
		\end{tabular}
\end{center}
\end{table}
Once trained, the neural network takes a total of $0.97$ hour for the first case during the optimization.  
A computationally efficient method described in De et al.~\cite{de2018computationally} can be used to solve \eqref{eq:nvie} but it will still take 1.33 days (projected) for function evaluations. On the other hand, a standard nonlinear solver, such as SciPy's \texttt{RK45}, takes $10.02$ s for one evaluation and a projected 5.2 days for $4.4\times10^4$ evaluations of the dynamical system during the optimization, which is impractical. 
For the second case, the proposed approach takes $1.09$ hours for optimization compared to 1.45 days (projected) for $4.8\times10^4$ evaluations of the dynamical system using the approach in De et al.~\cite{de2018computationally} and 5.67 days (projected) using SciPy's \texttt{RK45}. 
Hence, the proposed approach significantly reduces the computational cost of performing a robust optimization of nonlinear dynamical systems in this example. The reported computational times were obtained on a laptop equipped with an AMD Ryzen 7 8845HS processor (3.80 GHz), 16 GB of RAM, running 64-bit Windows 11.

\begin{figure}[!htb]
    \centering

    \begin{subfigure}[t]{0.48\textwidth}
        \centering
        \includegraphics[width=\textwidth]{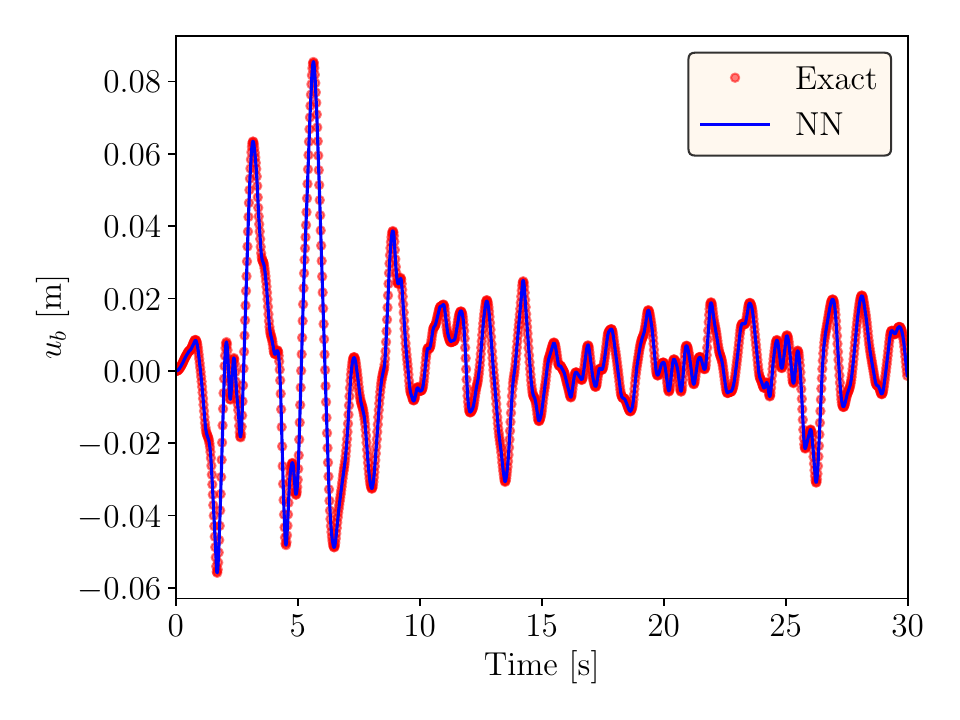}
        \caption{Base displacement, $u_b(t)$.}
        \label{fig:building_ub}
    \end{subfigure}
    \hfill
    \begin{subfigure}[t]{0.48\textwidth}
        \centering
        \includegraphics[width=\textwidth]{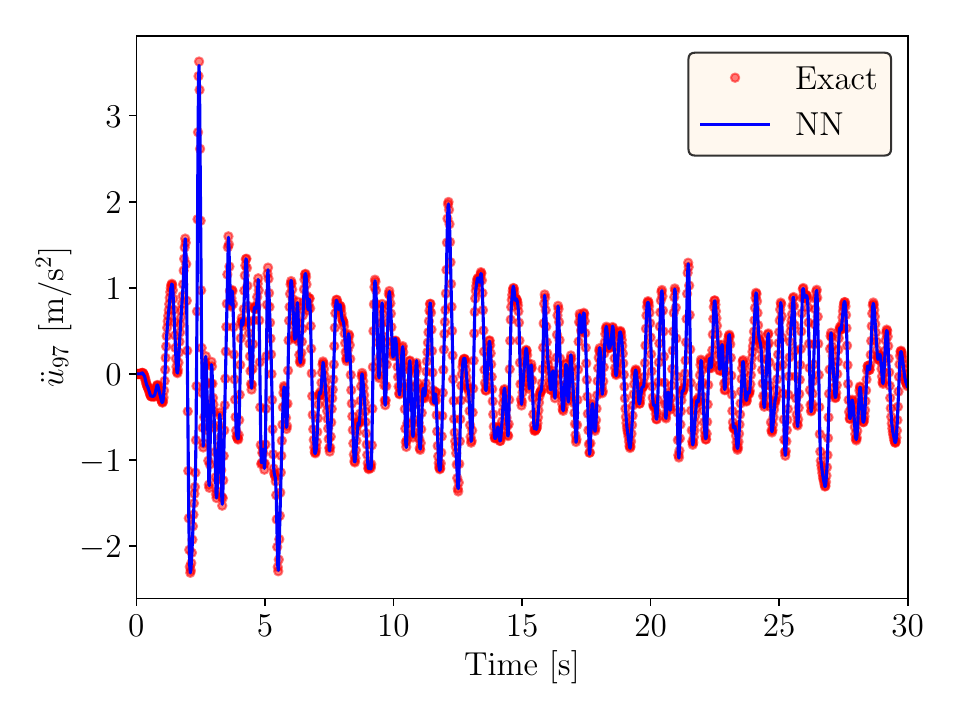}
        \caption{Horizontal roof acceleration, $\ddot{u}_{97}(t)$.}
        \label{fig:building_roof}
    \end{subfigure}

    \vspace{0.3cm}

    \begin{subfigure}[t]{0.48\textwidth}
        \centering
        \includegraphics[width=\textwidth]{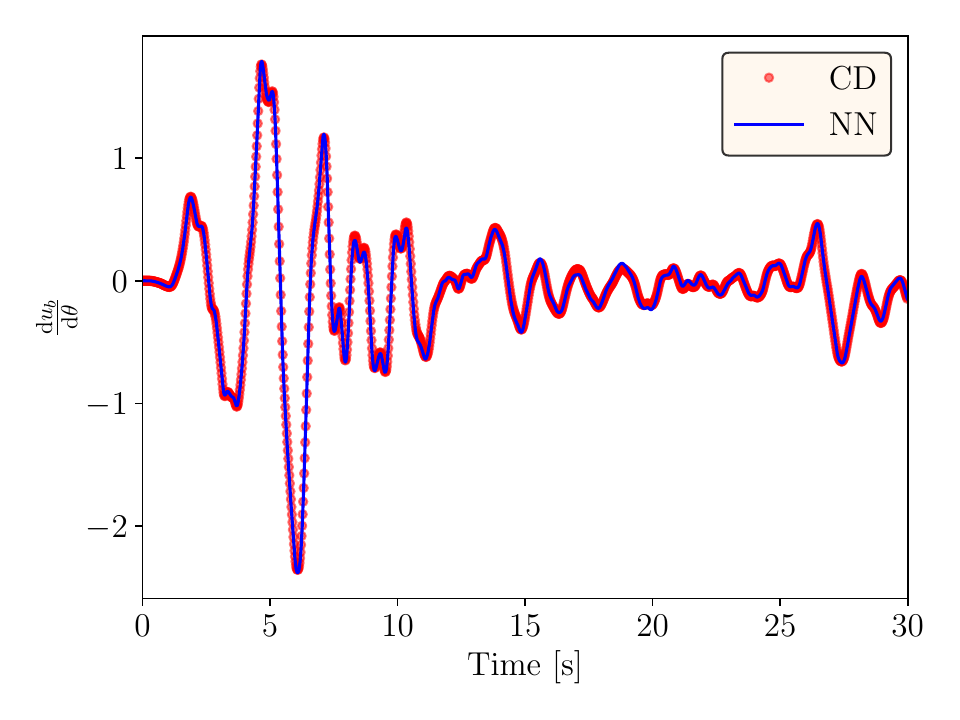}
        \caption{Base displacement sensitivity.}
        \label{fig:building_dub}
    \end{subfigure}
    \hfill
    \begin{subfigure}[t]{0.48\textwidth}
        \centering
        \includegraphics[width=\textwidth]{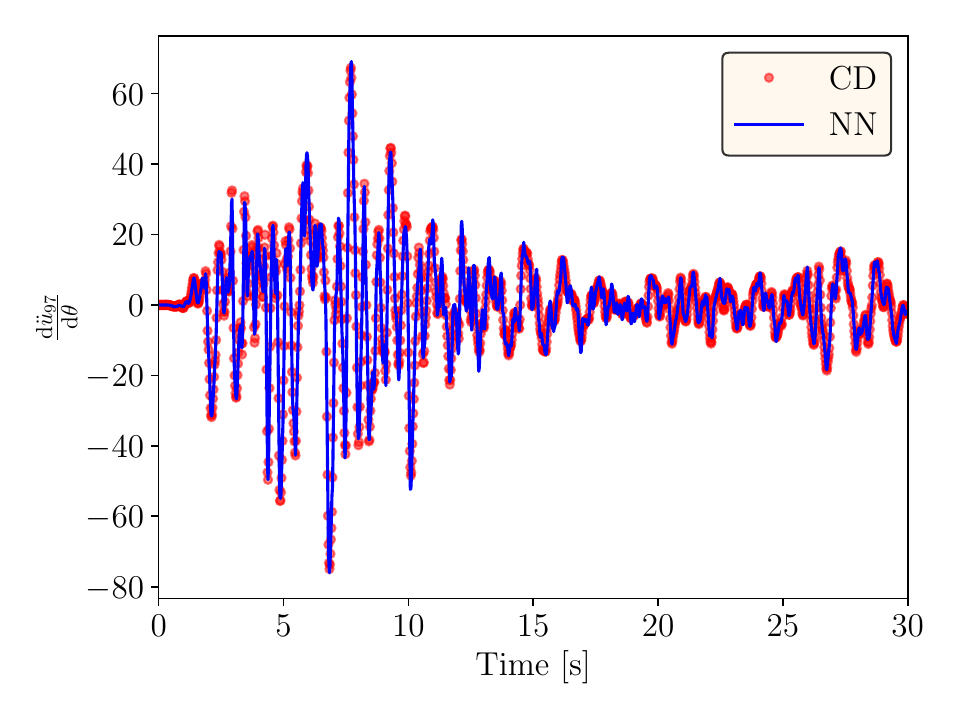}
        \caption{Horizontal roof acceleration sensitivity.}
        \label{fig:building_droof}
    \end{subfigure}

    \caption{Comparison of the proposed surrogate model
    and the high-fidelity nonlinear solution for Example II
    at the optimized isolator yield force. 
    The base displacement, horizontal roof acceleration, and their sensitivities with respect to design variable are evaluated for one realization of the uncertain parameters.
    }

    \label{fig:building_validation}
\end{figure}

The accuracy of the proposed surrogate model
is evaluated at the optimized design using 2,000
Monte Carlo samples. The statistical quantities obtained
from the proposed method are compared with those calculated
using the high-fidelity nonlinear model, 
in which the governing equations of motion, including the hysteretic behavior of the base isolation system defined in \eqref{eq:11story}, are solved directly without the neural network approximation, 
as reported in
Table~\ref{tab:opt_verification}. 
The mean components of the objective  obtained using the proposed and high-fidelity models are \(3.0649\) and \(3.0893\), respectively. The corresponding variance components are \(0.2831\) and \(0.2944\), while the mean base displacement constraint values are \(-0.0036\) and \(-0.0037\). 
These results indicate that the proposed surrogate model maintains reasonable accuracy in predicting the response quantities at the optimized
design under uncertainty. 
To further assess the surrogate accuracy at the response level, the time histories of the base displacement,
$u_b(t)$, and horizontal roof acceleration, $\ddot{u}_{97}(t)$, along with their sensitivities with respect to the design variable are compared against the high-fidelity nonlinear solution for one realization of the uncertain parameters. As illustrated in
Figure~\ref{fig:building_validation}, the proposed surrogate closely reproduces the base displacement and horizontal roof acceleration responses and its gradients.

\section{Conclusions} 

This paper developed a bi-fidelity method for the design of uncertain nonlinear dynamical systems with localized nonlinear behavior. The response from a nominal linear state-space model is combined with a neural network that approximates the nonlinear pseudoforce. The full response follows from the convolution of this pseudoforce with the impulse response operator of the linear system. As a result, the network does not have to reproduce the complete system dynamics. Additionally,  automatic differentiation provided derivatives of the pseudoforce with respect to the design variables, and the linear model propagated these derivatives to the response quantities used by the optimizer. The error analysis related the accuracy of the objective, constraints, and their gradients to errors in the pseudoforce approximation. This makes the pseudoforce and its derivative the relevant quantities to monitor during training and validation. 

In the solar heating example, the network modeled the nonlinear thermal contribution and its derivative with respect to the collector area. Each evaluation during optimization required only the neural network correction and the nominal linear model, rather than a full nonlinear simulation. As a result, the total computational time for optimization was reduced from hours to minutes. 
At the optimized design, verification using 2,000 Monte Carlo samples showed close agreement between the proposed method and the high-fidelity reference, with relative differences of approximately $0.005\%$ in the mean component of the objective and $0.1\%$ in the variance component of the objective. 
The base-isolated building provided a larger structural example with localized hysteresis. 
 In this example, the nonlinear uncertain pseudoforce was learned using a neural network, avoiding repeated solutions of the nonlinear differential equations during optimization. This reduced the total computational time for optimization from days to only a few hours. 
Similar to the first example, at the optimized design, verification using 2,000 Monte Carlo samples yielded to a relative difference of approximately $0.7\%$ in the mean component of the  objective. 

Although the proposed method substantially reduces the computational cost of optimization, it requires an offline training stage involving high-fidelity simulations. Furthermore, the feedforward neural network used in this study predicts the nonlinear pseudoforce at individual time steps. 
This limitation may affect its ability to capture complex nonlinear behavior with strong history dependence. Therefore, future work will investigate recurrent neural networks (RNNs) to better capture the temporal evolution of nonlinear dynamics and improve the accuracy of the predicted pseudoforce and its sensitivities. Extensions to reliability-based design optimization combined with adaptive sampling strategies will also be explored. 

\section*{Acknowledgments} 
The authors thank Dr. Mahmoud Kamalzare and Dr. Erik A. Johnson for kindly providing the finite element model of the 11-story base-isolated building used in Example II. 

\appendix

\section{Proof of Proposition 1} \label{sec:proof_p1} 


For each realization $\xii$, the Lipschitz continuity of $R$ gives
\begin{equation}
\left|R(\mathcal{Y})-R(\mathcal{Y}_{\mathrm{NN}})\right| \leq \Lambda_R\left\|\operatorname{vec}(\mathcal{Y})-\operatorname{vec}(\mathcal{Y}_{\mathrm{NN}})\right\|_2.
\end{equation}
Since $\mathcal{Y}$ and $\mathcal{Y}_{\mathrm{NN}}$ contain the outputs stacked over the time instances in $\mathcal{T}$,  
\begin{equation}
\left\|\operatorname{vec}(\mathcal{Y})-\operatorname{vec}(\mathcal{Y}_{\mathrm{NN}})\right\|_2 = \left(\sum_{t_k\in\mathcal{T}}\left\|\mathbf Y_k-\mathbf Y_{k,\mathrm{NN}}\right\|_2^2\right)^{1/2}.
\end{equation}
Using the inequality \(\|\mathbf v\|_2\leq\|\mathbf v\|_1\) for a vector with nonnegative entries,
\begin{equation}
\left\|\operatorname{vec}(\mathcal{Y})-\operatorname{vec}(\mathcal{Y}_{\mathrm{NN}})\right\|_2 \leq \sum_{t_k\in\mathcal{T}}\left\|\mathbf Y_k-\mathbf Y_{k,\mathrm{NN}}\right\|_2.
\end{equation}
For each $t_k\in\mathcal T$, using \eqref{eq:state_space_op}, $\Ym_k=\ym_k+\Delta t\sum_{i=0}^{k}\Cm\Hm_{L,k-i}\ppm_i$. Hence,  
\begin{equation} 
\begin{split}
    \left\|\Ym_k-\Ym_{k,\mathrm{NN}}\right\|_2&\leq\Delta t\sum_{i=0}^{k}\left\|\Cm\Hm_{L,k-i}\left(\ppm_i-\ppm_{i,\mathrm{NN}}\right)\right\|_2\\
&\leq\Gamma_C\Gamma_H\Delta t\sum_{i=0}^{k}\varepsilon_i.
\end{split} 
\end{equation} 
Combining the preceding inequalities gives the pointwise bound
\begin{equation}
\left|R(\mathcal{Y})-R(\mathcal{Y}_{\mathrm{NN}})\right| \leq \Lambda_R\Gamma_C\Gamma_H\Delta t\sum_{t_k\in\mathcal{T}}\sum_{i=0}^{k}\varepsilon_i.
\end{equation}
Taking the \(L^2(p_{\boldsymbol{\Xi}})\) norm of both sides gives
\begin{equation}
\left\|R(\mathcal{Y})-R(\mathcal{Y}_{\mathrm{NN}})\right\|_{L^2(p_{\boldsymbol{\Xi}})} \leq \Lambda_R\Gamma_C\Gamma_H\Delta t\left\|\sum_{t_k\in\mathcal{T}}\sum_{i=0}^{k}\varepsilon_i\right\|_{L^2(p_{\boldsymbol{\Xi}})}.
\end{equation}
Finally, Minkowski's inequality implies
\begin{equation}
\left\|\sum_{t_k\in\mathcal{T}}\sum_{i=0}^{k}\varepsilon_i\right\|_{L^2(p_{\boldsymbol{\Xi}})} \leq \sum_{t_k\in\mathcal{T}}\sum_{i=0}^{k}\left\|\varepsilon_i\right\|_{L^2(p_{\boldsymbol{\Xi}})}.
\end{equation}
Hence,
\begin{equation}
\left\|R(\mathcal{Y})-R(\mathcal{Y}_{\mathrm{NN}})\right\|_{L^2(p_{\boldsymbol{\Xi}})} \leq \Lambda_R\Gamma_C\Gamma_H\Delta t\sum_{t_k\in\mathcal{T}}\sum_{i=0}^{k}\left\|\varepsilon_i\right\|_{L^2(p_{\boldsymbol{\Xi}})},
\end{equation}
which proves the result. 
\hfill$\square$

\section{Proof of Corollary 1} \label{sec:proof_c1}
By the triangle inequality, 
\begin{equation}
\left|\mathcal R(\thetaa) - \mathcal R_{\mathrm{NN}}(\thetaa)\right| \leq \left|\Exp_{\xii}[Z] - \Exp_{\xii}[Z_{\mathrm{NN}}]\right| + \lambda_R \left|\operatorname{Var}_{\boldsymbol{\xi}}(Z) - \operatorname{Var}_{\xii}(Z_{\mathrm{NN}})\right|. 
\end{equation} 
The error in the mean is bounded by
\begin{equation}
\left|\Exp_{\xii}[Z] - \Exp_{\xii}[Z_{\mathrm{NN}}]\right| = \left|\Exp_{\xii} \left[Z-Z_{\mathrm{NN}}\right]\right|
\leq \left\|Z-Z_{\mathrm{NN}}\right\|_{L^2(p_{\boldsymbol{\Xi}})}.
\end{equation} 

Define the centered random variables $Z_c = Z-\Exp_{\xii}[Z]$, $Z_{\mathrm{NN},c} = Z_{\mathrm{NN}} - \Exp_{\xii}[Z_{\mathrm{NN}}]$. Then 
\begin{equation} 
\begin{split}
    \left|\operatorname{Var}_{\xii}(Z) - \operatorname{Var}_{\xii}(Z_{\mathrm{NN}})\right| 
    & = \left|\|Z_c\|_{L^2}^2 - \|Z_{\mathrm{NN},c}\|_{L^2}^2 \right|\\ 
    & \leq \left(\|Z_c\|_{L^2} + \|Z_{\mathrm{NN},c}\|_{L^2} \right) \left\|Z_c-Z_{\mathrm{NN},c}\right\|_{L^2}.
\end{split} 
\end{equation} 
Since centering is a contraction in $L^2$,
\begin{equation}
    \left\|Z_c-Z_{\mathrm{NN},c}\right\|_{L^2} \leq \left\|Z-Z_{\mathrm{NN}}\right\|_{L^2}.
\end{equation} 
Additionally, $\|Z_c\|_{L^2}=\sigma_Z$, $\|Z_{\mathrm{NN},c}\|_{L^2} = \sigma_{Z_{\mathrm{NN}}}$. 
Therefore, 
\begin{equation}
    \left|\operatorname{Var}_{\xii}(Z) - \operatorname{Var}_{\xii}(Z_{\mathrm{NN}})
\right| \leq \left(\sigma_Z+\sigma_{Z_{\mathrm{NN}}}\right) \left\|Z-Z_{\mathrm{NN}}\right\|_{L^2}. 
\end{equation}

Combining the mean and variance bounds gives 
\begin{equation}
\left|\mathcal R(\thetaa) - \mathcal R_{\mathrm{NN}}(\thetaa)\right| \leq \left[1+\lambda_R\left(\sigma_Z+\sigma_{Z_{\mathrm{NN}}}\right)\right]\left\|Z-Z_{\mathrm{NN}}\right\|_{L^2}.
\end{equation} 
Finally, Proposition~1 implies
\begin{equation}
    \left\|Z-Z_{\mathrm{NN}}\right\|_{L^2} \leq \delta_R(\thetaa), 
\end{equation}
which proves the result.
\hfill\(\square\)

\section{Proof of Proposition 2} \label{sec:proof_p2} 

Let $\am = \nabla_{\vey}R(\mathcal{Y})$, $\am_{\mathrm{NN}} = \nabla_{\vey} R(\mathcal{Y}_{\mathrm{NN}})$. 
Using the chain rule, $\nabla_{\thetaa}R(\mathcal{Y}) = \Sm_Y^{T}\am$, $\nabla_{\thetaa}R(\mathcal{Y}_{\mathrm{NN}}) = \Sm_{Y,\mathrm{NN}}^{T}\am_{\mathrm{NN}}$. 
Adding and subtracting $\Sm_{Y,\mathrm{NN}}^{T}\am$ gives

\begin{equation}
\nabla_{\thetaa}R(\mathcal{Y}) - \nabla_{\thetaa}R(\mathcal{Y}_{\mathrm{NN}}) = \left(\Sm_Y-\Sm_{Y,\mathrm{NN}}\right)^T\am + \Sm_{Y,\mathrm{NN}}^T \left(\am-\am_{\mathrm{NN}}\right).
\end{equation}
Therefore,
\begin{equation}
\left\|\nabla_{\thetaa}R(\mathcal Y) - \nabla_{\thetaa}R(\mathcal Y_{\mathrm{NN}})\right\|_{L^1}
\leq\left\|\|\Sm_Y-\Sm_{Y,\mathrm{NN}}\|_F\|\am\|_2\right\|_{L^1} + \left\|\|\Sm_{Y,\mathrm{NN}}\|_F\|\am-\am_{\mathrm{NN}}\|_2\right\|_{L^1}.
\end{equation}
Using the Cauchy-Schwarz inequality, 
\begin{equation}
\left\|\nabla_{\thetaa}R(\mathcal Y) - \nabla_{\thetaa}R(\mathcal Y_{\mathrm{NN}}) \right\|_{L^1}\leq \|\Sm_Y-\Sm_{Y,\mathrm{NN}}\|_{L^2} \|\am\|_{L^2} + \|\Sm_{Y,\mathrm{NN}}\|_{L^2} \|\am-\am_{\mathrm{NN}}\|_{L^2}.
\end{equation} 
The output sensitivity error satisfies
\begin{equation}
    \left\|\Sm_Y-\Sm_{Y,\mathrm{NN}}\right\|_{L^2}\leq\Delta t\Gamma_C\Gamma_H\sum_{t_k\in\mathcal T}\sum_{i=0}^{k}\left\|\varepsilon_i^S\right\|_{L^2}.
\end{equation}
The Lipschitz condition gives
\begin{equation}
    \left\|\am-\am_{\mathrm{NN}}\right\|_{L^2}\leq\Lambda_T\left\|\vey - \veyn\right\|_{L^2}, 
\end{equation}
and the convolution representation implies (see Appendix~\ref{sec:proof_p1})
\begin{equation}
    \left\|\vey - \veyn\right\|_{L^2}\leq\Gamma_C\Gamma_H\Delta t\sum_{t_k\in\mathcal T}\sum_{i=0}^{k}\left\|\varepsilon_i\right\|_{L^2}.
\end{equation}
Combining these bounds gives the result. 
\hfill$\square$  

\section{Proof of Corollary 2} \label{sec:proof_c2}
Proposition~2 gives
\begin{equation}
    \left\|\nabla_{\thetaa}R(\mathcal{Y}) - \nabla_{\thetaa}R(\mathcal{Y}_{\mathrm{NN}}) \right\|_{L^1(p_{\boldsymbol{\Xi}})} \leq B_1(\thetaa) \Delta t \Gamma_C\Gamma_H \sum_{t_k\in\mathcal T} \sum_{i=0}^{k} \left\|\varepsilon_i^S\right\|_{L^2(p_{\boldsymbol{\Xi}})} + \Lambda_T B_2(\thetaa) \Delta t \Gamma_C\Gamma_H \sum_{t_k\in\mathcal T} \sum_{i=0}^{k}
\left\|\varepsilon_i\right\|_{L^2(p_{\boldsymbol{\Xi}})}. 
\end{equation}
Using Assumption~1,
\begin{equation}
    \sum_{t_k\in\mathcal T}\sum_{i=0}^{k}\left\|\varepsilon_i^S \right\|_{L^2(p_{\boldsymbol{\Xi}})} \leq \Lambda_p \sum_{t_k\in\mathcal T} \sum_{i=0}^{k} \left\|\varepsilon_i\right\|_{L^2(p_{\boldsymbol{\Xi}})}. 
\end{equation}
Substituting this inequality into the bound from Proposition~2 yields
\begin{equation}
    \left\|\nabla_{\thetaa}R(\mathcal{Y}) - \nabla_{\thetaa}R(\mathcal{Y}_{\mathrm{NN}})\right\|_{L^1(p_{\boldsymbol{\Xi}})}
\leq \left[B_1(\thetaa)\Lambda_p+\Lambda_T B_2(\thetaa)\right]\Delta t \Gamma_C\Gamma_H \sum_{t_k\in\mathcal T} \sum_{i=0}^{k} \left\|\varepsilon_i\right\|_{L^2(p_{\boldsymbol{\Xi}})}. 
\end{equation}
This proves the result.
\hfill$\square$

\section{Proof of Corollary 3} \label{sec:proof_c3}
Under the stated regularity assumptions,
\begin{equation}
    \nabla_{\thetaa}\mathcal R = \Exp_{\xii}[\gm] + 2\lambda_R \Exp_{\xii}\left[\left(Z-\Exp_{\xii}[Z]\right)\left(\gm-\Exp_{\xii}[\gm]\right)\right].
\end{equation}
An analogous expression holds for $\nabla_{\thetaa}\mathcal R_{\mathrm{NN}}$.
The difference between the covariance terms can be written as 
\begin{equation}
    \operatorname{Cov}_{\xii}(Z,\gm) - \operatorname{Cov}_{\xii}(Z_{\mathrm{NN}},\gm_{\mathrm{NN}}) = \operatorname{Cov}_{\xii}(Z-Z_{\mathrm{NN}},\gm) + \operatorname{Cov}_{\xii}(Z_{\mathrm{NN}},\gm-\gm_{\mathrm{NN}}). 
\end{equation} 
The Cauchy-Schwarz inequality and the contraction property of centering in $L^2$ give
\begin{equation}
\begin{split}
    \left\|\operatorname{Cov}_{\xii}(Z-Z_{\mathrm{NN}},\gm)\right\|_2 &\leq \sigma_{\gm}\left\|Z-Z_{\mathrm{NN}}\right\|_{L^2}\\
    \left\|\operatorname{Cov}_{\xii}(Z_{\mathrm{NN}},\gm-\gm_{\mathrm{NN}})\right\|_2 & \leq \sigma_{Z_{\mathrm{NN}}}\left\|\gm-\gm_{\mathrm{NN}}\right\|_{L^2}
    \end{split}
\end{equation}
Furthermore, $\left\|\Exp_{\xii}[\gm] - \Exp_{\xii}[\gm_{\mathrm{NN}}]\right\|_2 \leq \left\|\gm-\gm_{\mathrm{NN}}\right\|_{L^2}$. 
By applying the argument in Appendix~\ref{sec:proof_p2} with the generalized H\"older inequality (instead of the Cauchy-Schwarz inequality) 
\begin{equation}
\|fg\|_{L^2} \leq \|f\|_{L^4}\|g\|_{L^4},
\end{equation}
the fourth moment assumptions in Corollary~3 give
\begin{equation}
\left\|\gm-\gm_{\mathrm{NN}}\right\|_{L^2}\leq\delta_R^S(\thetaa).
\end{equation} 
The result then follows by applying Proposition~\ref{prop1} to $\|Z-Z_{\mathrm{NN}}\|_{L^2}$.
\hfill$\square$ 
\bibliographystyle{unsrt} 
\bibliography{sample.bib}
\end{document}